\input amstex 
\documentstyle{amsppt} 
\loadbold
\let\bk\boldkey

\magnification=1200
\hsize=5.75truein
\vsize=8.75truein 
\hcorrection{.25truein}
\loadeusm \let\scr\eusm
\loadeurm 

\font\Rm=cmr12 
 
\font\Rrm=cmr17 
\font\fScr=eusb10 \define\Scr#1{\text{\fScr #1}} 
 
\define\Ind{\text{\rm Ind}}
\define\Aut#1#2{\text{\rm Aut}_{#1}(#2)}
\define\End#1#2{\text{\rm End}_{#1}(#2)}
\define\Hom#1#2#3{\text{\rm Hom}_{#1}({#2},{#3})} 
\define\GL#1#2{\roman{GL}_{#1}(#2)}
\define\M#1#2{\roman M_{#1}(#2)}
\define\Gal#1#2{\text{\rm Gal\hskip.5pt}(#1/#2)}
\define\Ao#1#2{\scr A^0_{#1}(#2)} 
\define\Go#1#2{\scr G_{#1}(#2)} 
\define\upr#1#2{{}^{#1\!}{#2}} 
\define\pre#1#2{{}_{#1\!}{#2}} 
\define\N#1#2{\text{\rm N}_{#1/#2}} 
\define\gwr#1{\frak G^\roman{wr}(#1)} 
\define\wgwr#1{\widehat{\frak G}^\roman{wr}(#1)} 
\define\sw{\text{\rm sw}} 

\let\ge\geqslant
\let\le\leqslant
\let\ups\upsilon
\let\vD\varDelta 
\let\ve\varepsilon 
\let\eps\epsilon 

\let\vF\varPhi 
\let\vG\varGamma

\let\vp\varpi

\let\vS\varSigma 
\let\vs\varsigma 

\let\vT\varTheta
\let\vU\varUpsilon 
\let\vX\varXi
\define\wP#1{\widehat{\scr P}_{#1}} 
\define\wW#1{\widehat{\scr W}_{#1}} 
\define\wss#1{\widehat{\scr W}^{\text{\rm ss}}_{#1}} 
\define\wwr#1{\widehat{\scr W}^{\text{\rm wr}}_{#1}} 
\define\wG#1{\widehat{\text{\rm GL}}_{#1}} 
\define\wA#1{{\widehat{\scr A\,}\!}_{#1}} 
\define\wR#1#2{\widehat{\scr R}_{#1}({#2})} 
\define\wRp#1#2{\widehat{\scr R}^+_{#1}({#2})} 
\document \baselineskip=14pt \parskip=4pt plus 1pt minus 1pt 
\topmatter \nologo \nopagenumbers
\title\nofrills \Rrm 
Higher ramification and the local Langlands correspondence 
\endtitle 
\rightheadtext{Higher ramification and Langlands correspondence} 
\author 
Colin J. Bushnell and Guy Henniart 
\endauthor 
\leftheadtext{C.J. Bushnell and G. Henniart}
\affil 
King's College London and Universit\'e de Paris-Sud 
\endaffil 
\address 
King's College London, Department of Mathematics, Strand, London WC2R 2LS, UK. 
\endaddress
\email 
colin.bushnell\@kcl.ac.uk 
\endemail
\address 
Laboratoire de Math\'ematiques d'Orsay, Univ\. Paris-Sud, CNRS, Universit\'e
Paris-Saclay, 91405 Orsay, France.
\endaddress 
\email 
Guy.Henniart\@math.u-psud.fr 
\endemail 
\date January 2016 \enddate 
\abstract 
Let $F$ be a non-Archimedean locally compact field. We show that the local Langlands correspondence over $F$ has a property generalizing the higher ramification theorem of local class field theory. If $\pi$ is an irreducible cuspidal representation of a general linear group $\roman{GL}_n(F)$ and $\sigma$ the corresponding irreducible representation of the Weil group $\Cal W_F$ of $F$, the restriction of $\sigma$ to a ramification subgroup of $\Cal W_F$ is determined by a truncation of the simple character $\theta_\pi$ contained in $\pi$, and conversely. Numerical aspects of the relation are governed by a Herbrand-like function $\Psi_\vT$ depending on the endo-class $\vT$ of $\theta_\pi$. We give a method for calculating $\Psi_\vT$ directly from $\vT$. Consequently, the ramification-theoretic structure of $\sigma$ can be predicted from the simple character $\theta_\pi$ alone. 
\endabstract 
\keywords Local Langlands correspondence, simple character, endo-class, ramification group, higher ramification theorem, Herbrand function  
\endkeywords 
\subjclassyear{2000}
\subjclass 22E50, 11S37, 11S15 \endsubjclass 
\endtopmatter 
\document \baselineskip=14pt \parskip=4pt plus 1pt minus 1pt
\subhead 
1 
\endsubhead 
We examine the local Langlands correspondence \cite{14}, \cite{17}, \cite{19}, \cite{21} for general linear groups over a non-Archimedean locally compact field $F$. We obtain striking new results connecting the fine structure of cuspidal representations of $\GL nF$, as in the classification scheme of \cite{9}, and the ramification-theoretic structure of Galois representations. 
\par 
Our main theorem generalizes the higher ramification theorem of local class field theory. It gives rise to a function analogous to the classical Herbrand function of a field extension. Our second theorem is an algorithm for calculating that function. Taken together, the results offer an unprecedented opportunity to transmit detailed structure across the correspondence, pointing a new direction for the subject. Here, we only indicate very first steps. 
\subhead 
2 
\endsubhead  
Let $\scr W_F$ be the Weil group of a separable algebraic closure $\bar F/F$. Let $\wW F$ be the set of equivalence classes of irreducible, smooth, complex representations of the locally profinite group $\scr W_F$. (From now on, when speaking of a representation of a locally profinite group, we will {\it always\/} assume it to be smooth and complex.) For each integer $n\ge 1$, let $\Ao nF$ be the set of equivalence classes of irreducible cuspidal representations of the general linear group $\GL nF$. To work in a dimension-free manner, we set $\wG F = \bigcup_{n\ge 1} \Ao nF$: given $\pi\in \wG F$, there is a unique integer $\roman{gr}(\pi) = m\ge 1$ such that $\pi \in \Ao mF$. 
\par
The local Langlands correspondence for $F$ provides a canonical bijection 
$$ 
\aligned 
\wG F &\longrightarrow \wW F, \\ \pi &\longmapsto \upr L\pi, 
\endaligned 
\tag 1
$$ 
such that $\dim \upr L\pi = \roman{gr}(\pi)$. The correspondence truly embodies a vast generalization of local class field theory. However, there is more to local class field theory than the existence of the Artin reciprocity map $\bk a_F: \scr W_F \to F^\times$. A mere existence statement falls short of revealing many useful properties and applications. So too for the Langlands correspondence: knowledge of its existence, or even a construction, does not automatically yield significant new insight. 
\subhead 
3 
\endsubhead 
An instance suggests itself. If $\eps$ is a real parameter, $\eps \ge 0$, let $\scr W_F^\eps$ be the corresponding ramification subgroup of $\scr W_F$ in the upper numbering convention of \cite{22}. In particular, $\scr W_F^0$ is the inertia subgroup $\scr I_F$ of $\scr W_F$. Let $\scr W_F^{\eps+}$ be the closure of the subgroup $\bigcup_{\delta>\eps} \scr W_F^\delta$. Thus $\scr W_F^{0+}$ is the wild inertia subgroup $\scr P_F$ of $\scr W_F$. 
\par 
The {\it first ramification theorem\/} of local class field theory asserts that $\bk a_F(\scr P_F)$ is the group $U^1_F = 1{+}\frak p_F$ of principal units in $F$. More generally, let $k\ge 1$ be an integer and write $U^k_F = 1{+}\frak p_F^k$. The {\it higher ramification theorem\/} asserts that $\bk a_F(\scr W_F^k) = U^k_F$ and $\bk a_F(\scr W_F^{k+}) = U^{1+k}_F$. It therefore yields an isomorphism between the group of characters of $U^k_F$ and the group of characters of $\scr W_F^k$ trivial on $\scr W_F^k \cap \scr W_F^\roman{der}$, where $\scr W_F^\roman{der} = \roman{Ker}\,\bk a_F$ is the (closed) derived subgroup of $\scr W_F$. Consequently, the fine structure of characters of $U^k_F$ is reflected in characters of $\scr W_F^k$. Of course, $\scr W_F^k$ admits characters that are not trivial on $\scr W_F^k \cap \scr W_F^\roman{der}$. 
\subhead 
4 
\endsubhead 
The Ramification Theorem of \cite{4} 8.2 Theorem, \cite{7} 6.1 provides a generalization of the first ramification theorem of local class field theory. It is written in terms of endo-classes of {\it simple characters\/} in $\GL nF$, in the sense of \cite{9} (and the Background notes below). Simple characters are very special characters of specific compact open subgroups of $\GL nF$, with a multitude of extraordinary properties. Not least is the ability to {\it transfer\/} simple characters between general linear groups of differing dimensions in a way that preserves relations of intertwining and conjugacy. This leads to the notion of {\it endo-equivalence\/} of simple characters, developed in \cite{2}. It provides an equivalence relation on the class of all simple characters in all general linear groups over $F$, the equivalence classes being called {\it endo-classes.} 
\par
A representation $\pi \in \Ao nF$ contains a unique conjugacy class of simple characters (Corollary 1 of \cite{6}). These necessarily lie in the same endo-class $\vT_\pi$. If $\sigma = \upr L\pi$, the Ramification Theorem asserts that the restriction $\sigma\,\big|\,{\scr P_F}$ of $\sigma$ to $\scr P_F$ depends only on $\vT_\pi$, and conversely. More precisely, if $\sigma\in \wW F$, then $\sigma\,\big|\,{\scr P_F}$ is a direct sum of irreducible representations of $\scr P_F$, all of which are $\scr W_F$-conjugate and occur with the same multiplicity. So, writing $\wP F$ for the set of equivalence classes of irreducible representations of the profinite group $\scr P_F$, the representation $\sigma$ yields a unique element of $\scr W_F\backslash \wP F$ that we choose to denote $[\sigma;0]^+$. On the other hand, let $\Scr E(F)$ be the set of endo-classes of simple characters over $F$. Given $\vT \in \Scr E(F)$, there exists $\pi \in \wG F$ so that $\vT = \vT_\pi$. If $\sigma = \upr L\pi$, the orbit $[\sigma;0]^+$ depends only on $\vT$ rather than the choice of $\pi$. We therefore denote it $\upr L\vT$. The map 
$$ 
\aligned 
\Scr E(F) &\longrightarrow \scr W_F\backslash \wP F, \\ \vT &\longmapsto \upr L\vT, 
\endaligned 
\tag 2 
$$ 
is then a bijection. Results of \cite{7}, \cite{3} show that the Langlands correspondence can, in essence, be re-constructed from the bijection (2) via an explicit process. 
\subhead 
5 
\endsubhead 
Our main result here shows how (2) may be refined into a {\it family\/} of bijections generalizing the higher ramification theorem of local class field theory. It is based on the fact that the Langlands correspondence preserves conductors of pairs. 
\par 
If $\sigma$ is a finite-dimensional, semisimple representation of $\scr W_F$, let $\sw(\sigma)$ be the Swan conductor of $\sigma$ and write $\vs(\sigma) = \sw(\sigma)/\dim \sigma$. For $\pi_1,\pi_2 \in \wG F$, let $\sw(\pi_1\times\pi_2)$ be the Swan conductor of the pair $(\pi_1,\pi_2)$. This is defined via the local constant $\ve(\pi_1\times\pi_2,s,\psi)$ of \cite{18}, \cite{23}. Setting 
$$ 
\vs(\pi_1\times\pi_2) = \frac{\sw(\pi_1\times\pi_2)}{\roman{gr}(\pi_1)\,\roman{gr}(\pi_2)}, 
$$ 
the correspondence (1) has the property 
$$ 
\vs(\pi_1\times \pi_2) = \vs(\upr L\pi_1\otimes \upr L\pi_2), \quad \pi_i \in \wG F. 
\tag 3
$$ 
\subhead 
6 
\endsubhead 
We exploit parallel structures carried by the sets $\Scr E(F)$ and $\scr W_F\backslash \wP F$. On the Galois side, one defines a pairing $\Delta$ on $\wW F$ by 
$$ 
\Delta(\sigma,\tau) = \roman{inf}\,\{\eps>0: \Hom{\scr W_F^\eps}\sigma\tau \neq 0\},\quad \sigma,\tau \in \wW F. 
$$ 
This is symmetric and satisfies an ultrametric inequality, but does not separate points. The value $\Delta(\sigma,\tau)$ depends only on the orbits $[\sigma;0]^+, [\tau;0]^+ \in \scr W_F\backslash \wP F$, so $\Delta$ induces a pairing, again denoted $\Delta$, on $\scr W_F\backslash \wP F$. The second version of $\Delta$ separates points and is an ultrametric on $\scr W_F\backslash \wP F$. The following result derives from \cite{15}. 
\proclaim{Proposition A} 
Let $\sigma \in \wW F$. There exists a unique continuous function $\vS_\sigma(x)$, $x\ge 0$, such that $\vs(\check\sigma \otimes \tau) = \vS_\sigma\big(\Delta(\sigma,\tau)\big)$, for all $\tau \in \wW F$. 
\endproclaim 
Here, $\check\sigma$ is the {\it contragredient\/} of $\sigma$. The {\it decomposition function\/} $\vS_\sigma(x)$ is given by a formula (3.1.2) expressing the way $\sigma$ decomposes when restricted to the ramification subgroups $\scr W_F^x$, $x>0$. Consequently, one needs detailed knowledge of the inner workings of $\sigma$ in order to write it down. It depends only on $[\sigma;0]^+$, so we sometimes write $\vS_\sigma = \vS_{[\sigma;0]^+}$. 
\par 
This material is covered in sections 1--3 and is mostly familiar, but we have taken care to ensure that the narrative is complete. A couple of deeper results have exact analogues on the GL-side. We have chosen to prove the GL-versions, in the appropriate place, and then deduce the Galois versions via the Langlands correspondence. 
\subhead 
7 
\endsubhead 
Rather more surprising is the existence of exact analogues on the $\roman{GL}$-side, developed in sections 4 and 5. That $\Scr E(F)$ carries a canonical ultrametric $(\vT,\vU) \mapsto \Bbb A(\vT,\vU)$ is already implicit in \cite{2}. It is given by an explicit formula (5.1.1) in terms of transfers of simple characters. However, the conductor formula of \cite{8} can be reformulated in terms of $\Bbb A$ to yield: 
\proclaim{Proposition B} 
Let $\vT \in \Scr E(F)$. There exists a unique continuous function $\vF_\vT(x)$, $x\ge 0$, such that $\vs(\check\pi \times \rho) = \vF_\vT\big(\Bbb A(\vT,\vT_\rho)\big)$, for any $\pi \in \wG F$ satisfying $\vT_\pi = \vT$ and any $\rho \in \wG F$. 
\endproclaim 
Again, $\check\pi$ is the contragredient of $\pi$. The {\it structure function\/} $\vF_\vT$ can be written down completely in terms of $\vT$ (4.4.1). Throughout sections 4 and 5, we have to pay attention to the behaviour of $\Bbb A$ and $\vF_\vT$ relative to tamely ramified extensions of the base field $F$. This prepares the way for later results. 
\par 
Propositions A and B are results of rather different kinds.  Proposition A, while not trivial, has no claim to great depth. Proposition B, on the other hand, emerges on combining two deep and highly developed theories, the complete account of the smooth dual of $\GL nF$ from \cite{9}, \cite{12}, \cite{13}, \cite{2}, and Shahidi's analysis of the Rankin-Selberg local constant in terms intertwining operators and Plancherel measure \cite{23}. This comparison is an instance of a common phenomenon: it is usually easier to access matters of depth via the $\roman{GL}$-side. 
\subhead 
8 
\endsubhead 
In section 6, we use (3) to combine the propositions and get the first of our main results. 
\proclaim{Higher Ramification Theorem} 
Let $\vT \in \Scr E(F)$. For $\eps > 0$, define $\delta>0$ by 
$$
\vF_\vT(\delta) = \vS_{\upr L\vT}(\eps). 
\tag 4 
$$ 
If $\vU \in \Scr E(F)$, then $\Bbb A(\vT,\vU) < \delta$ if and only if $\Delta(\upr L\vT,\upr L\vU) < \eps$. 
\endproclaim 
The result holds equally with non-strict inequalities. This form is easy to prove and contains everything of substance. However, working back through the definitions, one finds a more concrete version (6.5 Corollary). For representations $\sigma,\tau \in \wW F$, the condition $\Delta(\sigma,\tau) < \eps$ is equivalent to $\sigma$ and $\tau$ having a common irreducible component on restriction to $\scr W_F^\eps$. On the other side, take $\pi, \rho \in \wG F$. The condition $\Bbb A(\vT_\pi,\vT_\rho) < \delta$ is equivalent to $\pi$, $\rho$ each containing a representative of the same endo-class of truncated simple characters, in the more general sense of \cite{2}. The severity of the truncation is measured by $\delta$. The theorem thus implies a parametrization of conjugacy classes of representations of ramification groups by endo-classes of truncated simple characters, the Langlands correspondence inducing a bijection between the set of $\pi \in \wG F$ containing a given truncated endo-class and the set of $\sigma\in \wW F$ containing the corresponding representation of a ramification subgroup. 
\example{Example} 
Let $k\ge1$ be an integer and let $\phi$ be a character of $\scr W^k_F$ trivial on commutators: equivalently, $\phi = \tilde\phi\,\big|\,\scr W_F^k$, for some character $\tilde\phi$ of $\scr W_F$. Thus $\tilde\phi = \chi\circ \bk a_F$, for a character $\chi$ of $F^\times$. The restriction $\chi_1 = \chi\,\big|\,U^1_F$ is a simple character in $\GL1F$. The restriction $\chi_k = \chi\,\big|\,U^k_F$ is a truncation of $\chi_1$, and gives the endo-class corresponding to $\phi$ under the main theorem. 
\endexample 
\subhead 
9 
\endsubhead 
Our second main result concerns the change of scale $\eps \mapsto \delta$ in the Higher Ramification Theorem. Define a function $\Psi_\vT(x)$, $x\ge 0$, $\vT \in \Scr E(F)$, by $\Psi_\vT = \vF_\vT^{-1}\circ \vS_{\upr L\vT}$. Thus, in the theorem, $\delta = \Psi_\vT(\eps)$. The function $\Psi_\vT$ is continuous, positive, strictly increasing, piece-wise linear and smooth outside of a finite set. It plays a r\^ole analogous to the classical Herbrand functions, so we appropriate the name. Our second result, the {\it Interpolation Theorem\/} of section 7, gives a procedure for calculating $\Psi_\vT$ directly from $\vT$, without the recourse to the Langlands correspondence implicit in its definition. Since $\vT$ determines $\vF_\vT$ explicitly, the theorem yields the Galois-theoretic decomposition function $\vS_{\upr L\vT}$, {\it with no reference to Galois theory!} 
\par 
The Herbrand function $\Psi_\vT$ has simple behaviour relative to tamely ramified base field extension (7.1). Using this, we show that $\Psi_\vT$ can be calculated from the values of $\Bbb A(\vT, \chi\vT)$, as $\chi$  ranges over a certain set of characters of $F^\times$, along with the corresponding result relative to tame base field extensions. The final statement (7.5) is very simple, but the extraction of explicit formulas promises to be a challenging task. Here we examine only the easiest example of {\it essentially tame representations\/} (7.7). 
\par 
For our concluding section 8, we change to the Galois side to broach a related question: if we are given a decomposition function $\vS_\sigma$, what does it tell us about $\sigma$? We show that the first discontinuity of the derivative $\vS'_\sigma$ gives a canonical family of presentations of $\sigma$ as an induced representation in a manner respecting ramification structures. Recent work suggests this approach provides a useful complement to the Interpolation Theorem.  We end with a few specific examples. 
\subhead 
Background and notation 
\endsubhead
Throughout, $F$ is a non-Archimedean local field with finite residue field of characteristic $p$. The symbols $\frak o_F$, $\frak p_F$, $\Bbbk_F = \frak o_F/\frak p_F$, $U_F = \frak o_F^\times$, $U^k_F = 1{+}\frak p_F^k$, $k\ge 1$, and $\ups_F:F^\times \twoheadrightarrow \Bbb Z$ all have their customary meaning. 
\par 
Let $\bar F/F$ be a separable algebraic closure of $F$ and $\scr W_F = \scr W(\bar F/F)$ the Weil group of $\bar F/F$. Let $E/F$ be a finite separable extension of $F$. When working in the Galois-theoretic context, we generally assume $E$ to be a subfield of $\bar F$ and write $\scr W_E$ for the Weil group $\scr W(\bar F/E)$ of $\bar F/E$. We identify $\scr W_E$ with the open subgroup of $\scr W_F$ which fixes all elements of $E$ under the natural action of $\scr W_F$ on $\bar F$. 
\par 
We make extensive use of the theory of simple characters \cite{9}, along with endo-classes and tame lifting \cite{2}. An overview, containing what we need, can be found in \cite{1}: we give the barest summary here. 
\par 
Let $\frak a$ be a hereditary $\frak o_F$-order in $A = \End FV$, where $V$ is an $F$-vector space of finite dimension. We set $U_\frak a  = \frak a^\times$. If $\frak p$ is the Jacobson radical $\roman{rad}\,\frak a$ of $\frak a$, then $U^k_\frak a = 1{+}\frak p_\frak a^k$, $k\ge1$. We define the positive integer $e_\frak a$ by $\frak p_F\frak a = \frak p^{e_\frak a}$: this is the {\it $\frak o_F$-period of $\frak a$}. If $E/F$ is a subfield of $A$, we say $\frak a$ is {\it $E$-pure\/} if $x\frak ax^{-1} = \frak a$, for all $x\in E^\times$. 
\par 
Let $[\frak a,n,0,\beta]$ be a simple stratum in $A$ (\cite{9} (1.5.5)): in particular, the algebra $E = F[\beta]$ is a field and $\frak a$ is $E$-pure. As in \cite{9} 3.1, one attaches to this stratum an open subgroup $H^1(\beta,\frak a)$ of $U^1_\frak a$, and writes $H^k(\beta,\frak a) = H^1(\beta,\frak a)\cap U^k_\frak a$, $k\ge1$. 
\par 
Take a character $\psi_F$ of $F$ of level one (to use the terminology of \cite{9}). This means that $\psi_F$ is trivial on $\frak p_F$, but not trivial on $\frak o_F$. Following Chapter 3 of \cite{9}, one attaches to $[\frak a,n,0,\beta]$ and $\psi_F$ a specific non-empty, finite set $\scr C(\frak a,\beta,\psi_F)$ of characters of the compact group $H^1(\beta,\frak a)$. These are the {\it simple characters in\/} $\Aut FV$ defined by $[\frak a,n,0,\beta]$. The dependence on $\psi_F$ is rather trivial, so we usually regard it as permanently fixed, and omit it from the notation: thus $\scr C(\frak a,\beta,\psi_F) = \scr C(\frak a,\beta)$. 
\par 
In the same situation, let $m$ be an integer, $0\le m < n$. The symbol $\scr C(\frak a,m,\beta)$ means the set of characters of $H^{m+1}(\beta,\frak a)$ obtained by restricting the characters in $\scr C(\frak a,\beta)$: thus $\scr C(\frak a,0,\beta) = \scr C(\frak a,\beta)$. We refer to the elements of sets $\scr C(\frak a,m,\beta)$ as {\it truncated\/} simple characters. The general theory of endo-equivalence in \cite{2} applies equally to truncated simple characters. 
\head \Rm 
1. Ramification groups 
\endhead 
We start with a sequence of three sections on the Weil group and its representations. This one provides a brief {\it aide m\'emoire\/} for basic ramification theory and introduces some non-standard notation. 
\subhead 
1.1 
\endsubhead  
Let $\scr I_F$ be the inertia subgroup of $\scr W_F$ and $\scr P_F$ the wild inertia subgroup. Attached to a real number $\eps \ge -1$ is the ramification subgroup $\scr W^\eps_F$ of $\scr W_F$. We use the upper numbering convention of \cite{22} Chapitre IV, so that $\scr W_F^{-1} = \scr W_F$ and $\scr W_F^0 = \scr I_F$. This traditional notation is typographically inconvenient so, from now on, we use 
$$ 
\scr R_F(\eps) = \scr W_F^\eps, \quad \eps \ge0. 
\tag 1.1.1 
$$
The definition of the ramification sequence gives the semi-continuity property 
$$
\scr R_F(\eps) = \bigcap_{\delta<\eps} \scr R_F(\delta), \quad \eps>0. 
$$ 
One forms the subgroup $\bigcup_{\delta>\eps} \scr R_F(\delta)$ and its closure 
$$
\scr R_F^+(\eps) = \roman{cl}\big(\bigcup_{\delta>\eps} \scr R_F(\delta)\big) 
$$ 
in $\scr W_F$. This need not equal $\scr R_F(\eps)$: one says that $\eps$ is a {\it jump of $\bar F/F$} if $\scr R_F^+(\eps) \neq \scr R_F(\eps)$. In particular, 
$$ 
\scr R_F(0) = \scr I_F, \quad \scr R_F^+(0) = \scr P_F. 
$$
Each of the groups $\scr R_F(x)$, $\scr R_F^+(x)$, $x\ge0$, is profinite, closed and normal in $\scr W_F$. We summarize the main properties of the ramification groups, relative to finite quotients of $\scr W_F$, in the form we shall use them. We use the conventions of \cite{22} for numbering the ramification subgroups of a finite Galois group. 
\proclaim{Lemma} 
Let $x\ge 0$. Let $E/F$ be a finite Galois extension with $\vG = \Gal EF$. 
\roster 
\item 
The canonical image of $\scr R_F(x)$ in $\vG$ is the ramification group $\vG^x$. 
\item 
Suppose $x$ is not a jump in the ramification sequence for $E/F$, that is, $\vG^x= \vG^{x+\eps}$ for some $\eps >0$. The image of\/ $\scr R_F^+(x)$ in $\vG$ is then $\vG^x$. 
\item 
Suppose $x$ is a jump in the ramification sequence for $E/F$, that is, $\vG^x \neq \vG^{x+\eps}$, $\eps > 0$. 
\itemitem{\rm (a)} If $x$ is the largest jump for $E/F$, then the image of\/ $\scr R_F^+(x)$ in $\vG$ is trivial. 
\itemitem{\rm (b)} Otherwise, the image of\/ $\scr R_F^+(x)$ in $\vG$ is $\vG^y$, where $y$ is the least jump such that $y>x$. 
\endroster 
\endproclaim 
In the context of the lemma, it is often useful to have the notation $\vG^{x+} = \bigcup_{y>x} \vG^y$. Thus $x$ is a jump for $E/F$ if $\vG^x \neq \vG^{x+}$. In all cases, $\vG^{x+}$ is the image of $\scr R_F^+(x)$. 
\subhead 
1.2 
\endsubhead 
We make frequent use of the following facts. 
\proclaim{Lemma 1} 
If $K/F$ is a finite, tamely ramified field extension with $e = e(K|F)$, then  $\scr P_K = \scr P_F$ and 
$$ 
\alignedat3 
\scr R_F(x) &= \scr R_K(ex), &\quad x&> 0, \\ 
\scr R_F^+(x) &= \scr R_K^+(ex), &\quad x&\ge0. 
\endalignedat 
$$ 
\endproclaim 
\demo{Proof} 
This follows from the definition of the upper numbering of ramification groups. \qed 
\enddemo 
\proclaim{Lemma 2} 
If\/ $0<x\le y$, the commutator group $[\scr R_F(x), \scr R_F(y)]$ is contained in $\scr R_F^+(y)$. Moreover,  
$$ 
[\scr R_F^+(0),\scr R_F(x)] \subset \scr R_F^+(x), \quad x\ge0. 
$$
In particular, the group $\scr R_F(x)/\scr R_F^+(x)$, $x>0$, is central in $\scr R_F^+(0)/\scr R_F^+(x)$. 
\endproclaim 
\demo{Proof} 
The first assertion is implied by IV \S2 Proposition 10 of \cite{22}. The second then follows from the definition of $\scr R_F^+(x)$. \qed 
\enddemo 
\head \Rm 
2. Representations and ramification 
\endhead 
Let $\wW F$ be the set of isomorphism classes of irreducible representations of $\scr W_F$. Let $\wW F^\roman{ss}$ be the set of isomorphism classes of finite-dimensional semisimple representations of $\scr W_F$ ({\it cf\.} \cite{5} 28.7 Proposition). 
\par 
Let $\wR F\eps$ be the set of isomorphism classes of irreducible representations of the profinite group $\scr R_F(\eps)$, $\eps > 0$, and define $\wRp F\eps$, $\eps \ge 0$, analogously. The group $\scr W_F$ acts on both $\wR F\eps$ and $\wRp F\eps$ by conjugation. 
\par 
We investigate interactions between representations of $\scr W_F$ and the filtration by ramification groups. We identify the jumps in the ramification sequence for $\bar F/F$ and define a canonical pairing on $\wW F$. 
\subhead 
2.1 
\endsubhead 
We start at a general level. 
\proclaim{Proposition 1} 
An open subgroup $H$ of $\scr P_F$ contains $\scr R_F(\eps)$, for some $\eps>0$. Moreover, $H$ contains $\scr R_F(\eps')$, for some $\eps'<\eps$.  
\endproclaim 
\demo{Proof} 
The group $H$ is of the form $\scr P_F\cap \scr W_K$, for a finite extension $K/F$. Since $\scr R_F(\eps)$ is normal in $\scr W_F$, we may replace $H$ by the intersection of its $\scr W_F$-conjugates and assume $K/F$ is a Galois extension. If $\delta$ is the largest ramification jump for $K/F$, then any pair $\eps$, $\eps'$ satisfying $\delta <\eps'<\eps$ has the required property. \qed 
\enddemo 
\proclaim{Proposition 2} 
Let $\xi \in \wR F\eps$, $\eps > 0$. 
\roster 
\item 
The kernel of $\xi$ contains $\scr R_F(\delta)$, for some $\delta > \eps$. 
\item 
There exists $\sigma \in \wW F$ such that $\xi$ is equivalent to an irreducible component of the restriction $\sigma\,\big|\,{\scr R_F(\eps)}$. 
\endroster 
\endproclaim 
\demo{Proof} 
Since $\scr R_F(\eps)$ is profinite, the kernel of $\xi$ is open, hence of the form $H\cap \scr R_F(\eps)$ for an open subgroup $H$ of $\scr W_F$. Part (1) now follows from Proposition 1. 
\par 
If $E/F$ is a finite extension, set $\scr G_E  = \Gal{\bar F}E$. In part (2), it is enough to find an irreducible representation of $\scr G_F$ containing $\xi$ on restriction to $\scr R_F(\eps)$. We form the representation $I = \Ind_{\scr R_F(\eps)}^{\scr G_F}\,\xi$ of $\scr G_F$ smoothly induced from $\xi$. Thus $I$ is the union $\bigcup_{E/F} I^{\scr G_E}$ of its spaces of $\scr G_E$-fixed points, as $E/F$ ranges over finite Galois extensions contained in $\bar F$. The space $I^{\scr G_E}$ provides a representation of the finite group $\Gal EF$. Consequently, $I$ has an irreducible $\scr G_F$-subspace $\sigma$, and any such $\sigma$ has the desired property. \qed 
\enddemo 
\remark{Complement} 
Proposition 2 holds, with the same proof, on replacing $\wR F\eps$, $\eps >0$ with $\wRp F\eps$, $\eps \ge 0$. 
\endremark 
\remark{Apology} 
Proposition 2, applied to $\scr P_F = \scr R_F^+(0)$, replaces the incorrect proof of \cite{7} 1.2 Proposition. It also plugs a gap inadvertently left in the proof of \cite{4} 8.2 Theorem: we thank A\. Kilic for drawing our attention to the problem.  
\endremark 
\subhead 
2.2 
\endsubhead 
Let $\sigma \in \wW F$ and $\eps>0$. The restriction $\sigma\,\big|\,{\scr R_F(\eps)}$ of $\sigma$ to $\scr R_F(\eps)$ is semisimple. Its irreducible components are all $\scr W_F$-conjugate and occur with the same multiplicity. Thus $\sigma$ determines a unique conjugacy class $[\sigma;\eps] \in \scr W_F\backslash \wR F\eps$. Similarly, for $\eps \ge0$, $\sigma$ determines a unique conjugacy class $[\sigma;\eps]^+ \in \scr W_F\backslash \wRp F\eps$. 
\proclaim{Proposition} 
The orbit maps 
$$ 
\aligned 
\wW F &\longrightarrow \scr W_F\backslash \wR F\eps, \\ 
\sigma &\longmapsto [\sigma;\eps], 
\endaligned \qquad \eps > 0, 
$$ 
and 
$$ 
\aligned 
\wW F &\longrightarrow \scr W_F\backslash \wRp F\eps, \\ 
\sigma &\longmapsto [\sigma;\eps]^+, 
\endaligned \qquad \eps \ge 0, 
$$ 
are surjective. 
\endproclaim 
\demo{Proof} 
The assertion re-states 2.1 Proposition 2 and its complement. \qed 
\enddemo 
\subhead 
2.3 
\endsubhead 
Let $\sigma \in \wW F$. By 2.1 Proposition 1, $\roman{Ker}\,\sigma$ contains $\scr R_F(\eps)$ for $\eps$ sufficiently large. One defines the {\it slope\/} $\roman{sl}(\sigma)$ of $\sigma$ by 
$$ 
\roman{sl}(\sigma) = \roman{inf}\,\{\eps>0: \scr R_F(\eps) \subset \roman{Ker}\,\sigma\}. 
\tag 2.3.1 
$$ 
Thus $\roman{sl}(\sigma) = 0$ if and only if $\sigma$ is trivial on $\scr P_F$: one says that $\sigma$ is {\it tamely ramified.} 
\proclaim{Proposition} 
Let $\sigma \in \wW F$ and suppose that $\roman{sl}(\sigma) = s > 0$. The group $\scr R_F^+(s)$ is then contained in $\roman{Ker}\,\sigma$ while $\sigma\,\big|\,{\scr R_F(s)}$ is a direct sum of non-trivial characters of $\scr R_F(s)$. 
\endproclaim 
\demo{Proof} 
The first assertion follows from the definition of $\scr R_F^+(s)$. The group $\scr R_F(s)/\scr R_F^+(s)$ is abelian by 1.2 Lemma 2, so $\sigma\,\big|\,{\scr R_F(s)}$ is a direct sum of $\scr W_F$-conjugate characters. If these characters were trivial, $\scr R_F(s)$ would be contained in $\roman{Ker}\,\sigma$. Since $\roman{Ker}\,\sigma\cap \scr P_F$ is open in $\scr P_F$, it would contain $\scr R_F(t)$, for some $t<s$, by 2.1 Proposition 1, contrary to the definition of $s$. \qed 
\enddemo 
\proclaim{Corollary} 
\roster 
\item 
If $s>0$ is the slope of a representation $\sigma \in \wW F$, then $\scr R_F(s) \neq \scr R_F^+(s)$. In particular, $s$ is a jump in the ramification sequence for $\bar F/F$. 
\item 
If $s>0$ is a jump in the ramification sequence for $\bar F/F$, there exists $\sigma\in \wW F$ with slope $s$. 
\endroster 
\endproclaim 
\demo{Proof} 
Assertion (1) follows directly from the proposition. The profinite group $\scr R_F(s)$ admits a non-trivial smooth character $\xi$ which is trivial on the closed subgroup $\scr R_F^+(s)$. Assertion (2) is therefore given by 2.2 Proposition. \qed 
\enddemo 
\subhead 
2.4 
\endsubhead 
We can now identify the jumps in the ramification sequence, knowing that they all arise as slopes of irreducible representations. 
\par 
If $\rho\in \wW F^\roman{ss}$, let $\sw(\rho)$ denote the exponential {\it Swan conductor\/} of $\rho$. Thus $\sw(\rho)$ is a non-negative integer and, if we write $\rho = \bigoplus_{i=1}^r \tau_i$, with $\tau_i \in \wW F$, then $\sw(\rho) = \sum_{i=1}^r \sw(\tau_i)$. 
\proclaim{Basic connection} 
If $\sigma\in \wW F$, then $\roman{sl}(\sigma) = \sw(\sigma)/\dim\sigma$. In particular, $\roman{sl}(\sigma) \in \Bbb Q$. 
\endproclaim 
\demo{Proof} 
See Th\'eor\`eme 3.5 of \cite{16}. \qed 
\enddemo  
We complete the argument with a sharp result, which seems to lie rather deep. 
\proclaim{Proposition} 
Let $x>0$, $x\in \Bbb Q$, and write $x=a/b$, for relatively prime, positive integers $a$, $b$. There exists $\sigma\in \wW F$ such that $\sw(\sigma) = a$ and $\dim\sigma = b$. 
\endproclaim 
We defer the proof to 6.3 below. We deduce: 
\proclaim{Corollary} 
If $x\in \Bbb R$, $x > 0$, then $\scr R_F(x) \neq \scr R_F^+(x)$ if and only if $x\in \Bbb Q$. 
\endproclaim 
\subhead 
2.5 
\endsubhead 
The orbit maps $\wW F\to \scr W_F\backslash \wR F\eps$ and $\wW F\to \scr W_F\backslash \wRp F\eps$ of 2.2 factor through the orbit map $\wW F \to \scr W_F\backslash \wP F$. We use the same notation for the implied maps $\scr W_F\backslash \wP F\to \scr W_F\backslash \wR F\eps$ and $\scr W_F\backslash \wP F\to \scr W_F\backslash \wRp F\eps$. We set 
$$ 
\Delta(\xi,\zeta) = \roman{inf}\,\big\{\eps>0: [\xi;\eps] = [\zeta;\eps]\big\}, \quad \xi,\zeta\in \scr W_F\backslash \wP F. 
\tag 2.5.1 
$$ 
The pairing $\Delta$ is clearly symmetric: $\Delta(\xi,\zeta) = \Delta(\zeta,\xi)$. Its values are non-negative rational numbers (2.4 Corollary). Further, if $\delta = \Delta(\xi,\zeta) > 0$, then $[\xi;\delta]^+ = [\zeta;\delta]^+$ while $[\xi;\delta] \neq [\zeta;\delta]$. 
\proclaim{Proposition} 
\roster 
\item If $\xi,\zeta\in \scr W_F\backslash \wP F$, then $\Delta(\xi,\zeta) = 0$ if and only if $\xi = \zeta$. 
\item If $\xi,\zeta,\psi\in \scr W_F\backslash \wP F$, then 
$$ 
\Delta(\xi,\zeta) \le \roman{max}\,\big\{\Delta(\xi,\psi),\Delta(\psi,\zeta)\big\}. 
\tag 2.5.2 
$$ 
\endroster 
The pairing $\Delta$ is an ultrametric on $\scr W_F\backslash \wP F$. 
\endproclaim 
\demo{Proof} 
In part (1), one implication is trivial, so suppose $\xi\neq \zeta$. As in 2.2, there exists an irreducible representation $\tilde\xi$ of $\Gal{\bar F}F$ containing $\xi$ on restriction to $\scr P_F$. Choose $\tilde\zeta$ similarly. There exists a finite Galois extension $K/F$ such that both $\tilde\xi$ and $\tilde\zeta$ are inflated from representations of $\Gal KF$. The extension $K/F$ is not tamely ramified, so it has a least positive ramification jump $\phi$. If $0<\eps<\phi$, then $[\xi;\eps] \neq [\zeta;\eps]$, and this implies $\Delta(\xi,\zeta) \ge \eps >0$. Part (2) follows directly from the definition. \qed 
\enddemo 
It is often more convenient to view $\Delta$ as a pairing on $\wW F$, setting 
$$ 
\Delta(\sigma,\tau) = \Delta\big([\sigma;0]^+,[\tau;0]^+\big), \quad \sigma,\tau \in \wW F. 
$$ 
This version is again symmetric and has the ultrametric property (2.5.2), but it does not separate points. In this form, 
$$ 
\Delta(\sigma,\tau) = \roman{inf}\,\big\{\eps>0: \Hom{\scr R_F(\eps)}\sigma\tau \neq 0\big\}. 
\tag 2.5.3 
$$ 
\subhead 
2.6 
\endsubhead 
We consider the behaviour of $\Delta$ under tamely ramified base field extension. Temporarily write $\Delta_F$ for the pairing (2.5.1). Let $K/F$ be a finite tame extension, and let $\Delta_K$ be the analogue of $\Delta_F$ relative to the base field $K$. Thus $\scr P_K = \scr P_F$ and $\scr R_K(\eps) = \scr R_F(\eps/e)$, where $e = e(K|F)$ (1.2 Lemma 1). Consequently: 
\proclaim{Proposition} 
If $\xi, \zeta \in \wP F = \wP K$, then 
$$ 
e\,\Delta_F(\scr W_F\xi,\scr W_F\zeta) = \roman{min}\,\big\{\Delta_K(\scr W_K\xi, \scr W_Kg\zeta): g\in \scr W_K\backslash \scr W_F\big\}. 
$$ 
\endproclaim 
\head \Rm 
3. Ultrametric and conductors
\endhead 
We link the ultrametric $\Delta$ on $\scr W_F\backslash \wP F$ to conductors of tensor products of representations of $\scr W_F$. The basic ideas come from Heiermann's note \cite{15}. 
\subhead 
3.1  
\endsubhead 
Set 
$$ 
\vs(\sigma) = \sw(\sigma)/\dim\sigma, \quad \sigma\in \wW F^\roman{ss}. 
$$ 
If $\sigma \in \wW F$, this reduces to $\vs(\sigma) = \roman{sl}(\sigma)$, as in 2.4. 
\par 
Let $\sigma \in \wW F$, say $\sigma:\scr W_F \to \Aut{\Bbb C}V$, for a finite dimensional complex vector space $V$. The semisimple representation $\check\sigma\otimes \sigma$ thus acts on the space $X_\sigma = \check V\otimes V$. Write 
$$ 
X_\sigma(\delta) = X_\sigma^{\scr R_F^+(\delta)} = \bigcap_{\eps>\delta} X_\sigma^{\scr R_F(\eps)}, \quad \delta \ge 0, 
$$ 
for the space of $\scr R_F^+(\delta)$-fixed points in $X_\sigma$. Let $X'_\sigma(\delta)$ be the unique $\scr R_F^+(\delta)$-complement of $X_\sigma(\delta)$ in $X_\sigma$. Since $\scr R^+_F(\delta)$ is a normal subgroup of $\scr W_F$, the spaces $X_\sigma(\delta)$, $X'_\sigma(\delta)$ are $\scr W_F$-stable and provide semisimple representations of $\scr W_F$. 
\proclaim{Lemma} 
Let $\sigma,\tau \in \wW F$. If $\delta = \Delta(\sigma,\tau)$,  then 
$$
\vs(\check\sigma\otimes \tau) = (\dim\sigma)^{-2}\big(\delta\dim X_\sigma(\delta) + \sw\,X'_\sigma(\delta)\big). 
\tag 3.1.1 
$$ 
\endproclaim 
This formulation is to be found on p\. 572 of \cite{15} {\it cf\.} (3.1.3) below. We need a slightly different emphasis. 
\proclaim{Proposition} 
For $\sigma\in \wW F$ and $\delta\ge0$, define 
$$
\vS_\sigma(\delta) = (\dim\sigma)^{-2}\big(\delta\dim X_\sigma(\delta) + \sw\,X'_\sigma(\delta)\big). 
\tag 3.1.2
$$ 
The function $\vS_\sigma$ is continuous, strictly increasing, piecewise linear and convex. Its derivative $\vS'_\sigma$ is continuous outside of a finite set. 
\endproclaim 
\demo{Proof} 
Write $\check\sigma\otimes\sigma = \sum_i\psi_i$, where the $\psi_i$ are irreducible. The basic connection recalled in 2.4 implies that  
$$
\delta\dim X_\sigma(\delta) + \sw\,X'_\sigma(\delta) = \sum_i \roman{max}\,\{\delta\dim\psi_i, \sw\,\psi_i\}, 
\tag 3.1.3  
$$ 
and also that each term in the sum is a continuous, non-decreasing, function. One factor $\psi_i$ is the trivial representation, and that contributes a strictly increasing term. All assertions are now immediate. \qed 
\enddemo 
Comparing with (3.1.1), we have 
$$ 
\vS_\sigma\big(\Delta(\sigma,\tau)\big) = \vs(\check\sigma\otimes\tau), \quad \sigma,\tau \in \wW F. 
\tag 3.1.4 
$$ 
There is a consequence, useful in more general applications, although we do not need it here. 
\proclaim{Corollary} 
The pairing $(\sigma,\tau) \mapsto \vs(\check\sigma\otimes \tau)$ of \rom{(3.1.4)} satisfies the ultrametric inequality 
$$ 
\vs(\check\sigma_1\otimes\sigma_2) \le \roman{max}\,\{\vs(\check\sigma_1 \otimes \sigma_3),\vs(\check\sigma_3\otimes \sigma_2)\}, \quad \sigma_i \in \wW F. 
$$
\endproclaim 
\demo{Proof} 
The proof is identical to that of 5.4 Corollary below. We choose to give the details there. \qed 
\enddemo 
\definition{Notation} 
The function $\vS_\sigma$, as defined in (3.1.2), depends only on the class $[\sigma;0]^+ \in \scr W_F\backslash \wP F$. It is sometimes necessary to reflect this via the notation 
$$ 
\vS_\sigma(x) = \vS_{[\sigma;0]^+}(x). 
\tag 3.1.5 
$$ 
\enddefinition 
\subhead 
3.2 
\endsubhead 
Say that $\sigma\in \wW F$ is {\it totally wild\/} if the restriction $\sigma\,\big|\,{\scr P_F}$ is irreducible. If $\sigma$ is such a representation, and if $K/F$ is a finite, tamely ramified field extension, the restriction $\sigma^K = \sigma\,\big|\,{\scr W_K}$ is irreducible and totally wild. 
\proclaim{Proposition} 
If $\sigma \in \wW F$ is totally wild and $K/F$ is a finite tame extension then 
$$ 
\vS_\sigma(x) = e^{-1}\vS_{\sigma^K}(ex), \quad x\ge 0, 
$$ 
where $e = e(K|F)$. 
\endproclaim 
\demo{Proof} 
This follows from 1.2 Lemma 1. \qed 
\enddemo 
\subhead 
3.3 
\endsubhead 
The function $\vS_\sigma$ has a strong uniqueness property, although we don't need it at this stage. 
\proclaim{Proposition} 
The function $\vS_\sigma$, defined by \rom{(3.1.2)}, is the unique continuous function satisfying \rom{(3.1.4).} 
\endproclaim 
The proposition is an immediate consequence of the following, proved in 6.3 below.  
\proclaim{Density Lemma} 
Let $\sigma\in \wW F$. The set $\{\Delta(\sigma,\tau):\tau\in \wW F\}$ is dense on the half-line $x\ge 0$, $x\in \Bbb R$. 
\endproclaim 
\subhead 
3.4 
\endsubhead 
We record a continuity property of the function $\sigma \mapsto\vS_\sigma$.
\proclaim{Proposition} 
If $\sigma,\tau \in \wW F$ and $\delta = \Delta(\sigma,\tau)$, then $\vS_\sigma(x) = \vS_\tau(x)$, $x \ge \delta$. 
\endproclaim 
\demo{Proof} 
By definition, the condition $\delta = \Delta(\sigma,\tau)$ is equivalent to $[\sigma;\eps] = [\tau;\eps]$ for all $\eps >\delta$. If $\dim \sigma = a$, $\dim\tau = b$, this condition is equivalent to $b\sigma\,\big|\,{\scr R_F(\eps)} \cong a\tau\,\big|\,{\scr R_F(\eps)}$, for $\eps > \delta$. Comparing first trivial components and then non-trivial ones, we get 
$$
b^2X_\sigma(\eps) \cong a^2X_\tau(\eps), \quad b^2X'_\sigma(\eps) \cong a^2X'_\tau(\eps). 
$$
The assertion now follows from the definition (3.1.2). \qed 
\enddemo 
\head \Rm
4. Invariants of simple characters 
\endhead 
We pass to the $\roman{GL}$-side. In this section, we recall and develop some features of the theory of simple characters using mainly \cite{9} and \cite{2}. 
\subhead 
4.1 
\endsubhead 
We start with a detail from \cite{8} 6.4. Let $E/F$ be a finite field extension and let $A = \End FE$. Let $\frak a$ be the unique $E$-pure hereditary $\frak o_F$-order in $A$. Let $\beta \in E^\times$ satisfy $E = F[\beta]$ and $m = -\ups_E(\beta) > 0$. We assume that the quadruple $[\frak a,m,0,\beta]$ is a {\it simple stratum.} 
\par 
Let $a_\beta$ denote the adjoint map $A\to A$, $x\mapsto \beta x{-}x\beta$, and $s_{E/F}:A\to E$ a tame corestriction relative to $E/F$, \cite{9} 1.3. The sequence 
$$ 
0 \to E \longrightarrow A @>{\ a_\beta\ }>> A @>{\ s_{E/F}\ }>> E \to 0 
$$ 
is then exact. There exist $\frak o_F$-lattices $\frak l$, $\frak l'$ in $E$ and $\frak m$, $\frak m'$ in $A$ such that the sequence 
$$ 
0 \to \frak l \longrightarrow \frak m @>{\ a_\beta\ }>> \frak m' @>{\ s_{E/F}\ }>> \frak l' \to 0  
$$ 
is exact. For Haar measures $\mu_E$, $\mu_A$ on $E$ and $A$ respectively, the quantity 
$$
C(\beta) = \frac{\mu_E(\frak l)\,\mu_A(\frak m')}{\mu_E(\frak l')\,\mu_A(\frak m)} 
$$ 
is independent of all these choices. If $q = |\Bbbk_F|$, there is an integer $\frak c(\beta)$ such that 
$$ 
C(\beta) = q^{\frak c(\beta)}. 
\tag 4.1.1 
$$ 
As an example, consider the case where $\beta$ is {\it minimal over $F$.} In concrete terms, this means that $m = -\ups_E(\beta)$ is relatively prime to $e = e(E|F)$ and, for a prime element $\vp$ of $F$, the coset $\beta^e\vp^m{+}\frak p_E$ generates the residue field extension $\Bbbk_E/\Bbbk_F$. 
\proclaim{Proposition} 
Set $e = e(E|F)$, $f=f(E|F)$. If $\beta$ is minimal over $F$, then $\frak c(\beta) = mf(ef{-}1)$. 
\endproclaim 
\demo{Proof} 
In this situation, the sequence 
$$
0\to \frak o_E \longrightarrow \frak a @>{\ a_\beta\ }>> \beta\frak a @>{\ s_{E/F}\ }>> \beta\frak o_E \to 0 
$$ 
is exact \cite{9} (1.4.15). The result then follows from a short calculation. \qed 
\enddemo 
\subhead 
4.2 
\endsubhead 
Let $\Scr E(F)$ be the set of endo-classes of simple characters over $F$, including the trivial element $\bk 0$. Let $\vT \in \Scr E(F)$, $\vT \neq \bk 0$. There is a finite-dimensional $F$-vector space $V$, a simple stratum $[\frak a,m,0,\beta]$ in $\End FV$ and a simple character $\theta\in \scr C(\frak a,\beta,\psi_F)$ such that $\theta$ has endo-class $\vT$. One says that $\theta$ is a {\it realization of $\vT$} (on $\frak a$, relative to $\psi_F$). Let $e_\frak a$ be the $\frak o_F$-period of $\frak a$. We recall \cite{2} (8.11) that the quantities 
$$ 
\alignedat3 
\deg\vT &= [F[\beta]{:}F], &\quad m_\vT &= m/e_\frak a, \\ 
e_\vT &= e(F[\beta]|F), &\quad f_\vT &= f(F[\beta]|F) 
\endalignedat 
\tag 4.2.1 
$$ 
depend only on $\vT$ and not on the choices of $\theta$, $\frak a$, $\psi_F$ or $\beta$. The same applies to 
$$ 
k_0(\vT) = k_0(\beta,\frak a)/e_\frak a, 
\tag 4.2.2 
$$ 
where $k_0(\beta,\frak a)$ is the ``critical exponent'' of \cite{9} (1.4.5). Recall that $k_0(\vT) = -\infty$ when $\deg\vT = 1$. Otherwise, $k_0(\vT)$ is a negative rational number satisfying $-k_0(\vT) \le m_\vT$. 
\par 
If $\frak a$ is a hereditary order in $A= \End FV$, the realization of the trivial element $\bk 0$ of $\Scr E(F)$ on $\frak a$ is the trivial character of the group $U^1_\frak a = 1{+}\frak p$, where $\frak p = \roman{rad}\,\frak a$. We set 
$$ 
\gathered 
\deg\bk 0 = e_{\bk 0} = f_{\bk 0} = 1, \\
m_{\bk 0} = 0. 
\endgathered 
\tag 4.2.3 
$$ 
The following observation will be useful later. 
\proclaim{Proposition} 
Let $x$ be a positive rational number, say $x = a/b$, for relatively prime positive integers $a$, $b$. There exists $\vT \in \Scr E(F)$ such that $m_\vT = x$ and $e_\vT = \deg\vT = b$. 
\endproclaim 
\demo{Proof} 
Let $E/F$ be a totally ramified field extension of degree $b$ and choose $\alpha\in E$ of valuation $-a$. The element $\alpha$ is then minimal over $F$. If $\frak a$ is the unique $E$-pure hereditary $\frak o_F$-order in $\End FE$, the quadruple $[\frak a,a,0,\alpha]$ is a simple stratum. The endo-class $\vT$ of any $\theta\in \scr C(\frak a,\alpha)$ has the required properties. \qed 
\enddemo 
\subhead 
4.3 
\endsubhead 
Let $\vT \in \Scr E(F)$, $\vT \neq \bk0$. We attach to $\vT$ a finite set $\scr S_\vT$ of positive rational numbers, to be called the {\it normalized jumps\/} of $\vT$. We choose a realization $\theta\in \scr C(\frak a,\beta)$ of $\vT$, as in 4.2. We first attach to $[\frak a,m,0,\beta]$ a finite set $\roman S_{[\frak a,\beta]}$ of positive integers $t$, such that $-k_0(\beta,\frak a) \le t \le m$. 
\par 
We proceed by induction on the degree $[F[\beta]{:}F]$. If $\beta$ is minimal over $F$, in particular if $\beta\in F^\times$, we put $\roman S_{[\frak a,\beta]} = \{m\}$. If $\beta$ is not minimal over $F$, we set $r= -k_0(\beta,\frak a)$. Thus $0<r<m$. We choose a simple stratum $[\frak a,m,r,\gamma]$ equivalent to $[\frak a,m,r,\beta]$ \cite{9} (2.4.1). Thus $[\frak a,m,0,\gamma]$ is simple and $[F[\gamma]{:}F] < [F[\beta]{:}F]$. The set $\roman S_{[\frak a,\gamma]}$ has been defined inductively, and its least element is either $m$ or $-k_0(\gamma,\frak a)$. In either case, the least element is strictly greater than $r$. We define 
$$
\roman S_{[\frak a,\beta]} = \roman S_{[\frak a,\gamma]} \cup \{r\}. 
$$ 
\remark{Remark} 
If we have another simple stratum $[\frak a',m',0,\beta]$ in $\End F{V'}$, for some $V'$, then $\roman S_{[\frak a',\beta]} = \{xe_\frak a/e_{\frak a'}: x\in \roman S_{[\frak a,\beta]}\}$, as follows from \cite{9} (1.4.13). 
\endremark 
We define 
$$ 
\scr S_\vT = \{s/e_\frak a: s\in \roman S_{[\frak a,\beta]}\}. 
\tag 4.3.1 
$$ 
The least element of $\scr S_\vT$ is thus either $m_\vT$ or $-k_0(\vT)$. 
\proclaim{Lemma} 
The set $\scr S_\vT$ depends only on $\vT$, and not on any of the choices $\theta$, $\psi_F$, $[\frak a,m,0,\beta]$. 
\endproclaim 
\demo{Proof} 
This follows from \cite{9} (3.5.4). \qed 
\enddemo 
\definition{Definition} 
Let $x\in \Bbb R$, $x\ge 0$, $x\notin \scr S_\vT$. 
\roster 
\item 
If $x < \roman{min\,} \scr S_\vT$, set $\gamma_x = \beta$. 
\item If $x> m_\vT = \roman{max\,} \scr S_\vT$, set $\gamma_x = 0$. 
\item 
Otherwise, let $y = t/e_\frak a$ be the least element of $\scr S_\vT$ such that $y > x$, and let $[\frak a,m,t,\gamma_x]$ be a simple stratum equivalent to $[\frak a,m,t,\beta]$. 
\endroster 
Set $E_x = F[\gamma_x]$ and define 
$$
\aligned 
d_\vT(x) &= [E_x{:}F], \\ e_\vT(x) &= e(E_x|F), \\ \frak c_\vT(x) &= \frak c(\gamma_x). 
\endaligned 
\tag 4.3.2 
$$ 
\enddefinition 
\proclaim{Proposition} 
The quantities \rom{(4.3.2)} depend only on $x$ and $\vT$. If $y_1<y_2$ are successive elements of the set $\{0,\infty\}\cup \scr S_\vT$, the functions \rom{(4.3.2)} are constant in the region $y_1<x<y_2$. 
\endproclaim 
\demo{Proof} 
This follows, via an inductive argument, from the properties recalled in 4.2. \qed 
\enddemo 
\remark{\bf Observation} 
The proposition notwithstanding, all the invariants (4.3.2) of $\vT$ are defined purely in terms of an element $\beta$ giving rise to a realization of $\vT$. 
\endremark 
\subhead 
4.4 
\endsubhead 
Let $\vT\in \Scr E(F)$, $\vT \neq \bk0$, be as 4.3. We define a function $\vF_\vT(x)$, $x\ge0$. To start with, assume $x \notin \scr S_\vT$ and set  
$$ 
\vF_\vT(x) = \left\{\, \alignedat3 &x, &\qquad  &x > m_\vT, \\ 
&\frac{\frak c_\vT(x)}{d_\vT(x)^2} + \frac x{d_\vT(x)}, &\qquad &0 < x < m_\vT. 
\endalignedat \right. 
\tag 4.4.1 
$$ 
\proclaim{Proposition} 
\roster 
\item 
The function $\vF_\vT(x)$ of \rom{(4.4.1)} extends uniquely to a continuous function on the half-line $x\ge0$, that is, 
$$ 
\underset{x\to y-}\to{\roman{lim}}\,\vF_\vT(x) = \underset{x\to y+}\to{\roman{lim}}\,\vF_\vT(x), \quad y\in \scr S_\vT. 
$$ 
\item 
The function $\vF_\vT$ is piecewise linear, convex and strictly increasing. 
\item 
If $x\notin \scr S_\vT$, then $\vF_\vT'(y) = d_\vT(x)^{-1}$, for $y$ ranging over some open neighbourhood of $x$. 
\item 
The discontinuities of the derivative $\vF'_\vT$ are the elements $x$ of\/ $\scr S_\vT$ \rom{except} when $E_{m_\vT-\eps} = F$ for some $\eps>0$. In that case, $\vF'_\vT$ is continuous at $m_\vT$. 
\endroster 
\endproclaim 
\demo{Proof} 
Assertion (1) is given by 4.1 Proposition above together with 3.1 Proposition of \cite{4}. Part (2) then follows from (4.4.1) and 4.3 Proposition. Part (3) follows from the definition and 4.3 Proposition, part (4) from the definition.  \qed 
\enddemo 
The trivial element $\bk0$ of $\Scr E(F)$ is dealt with via the explicit formula 
$$ 
\vF_{\bk0}(x) = x,\quad x\ge0. 
\tag 4.4.2 
$$ 
In all cases, we call $\vF_\vT$ the {\it structure function\/} of $\vT\in \Scr E(F)$. 
\proclaim{Complements} 
Let $[\frak a,m,0,\beta]$ be a simple stratum. 
\roster 
\item 
For $i=1,2$, let $\theta_i\in \scr C(\frak a,\beta)$. If $\vT_i$ is the endo-class of $\theta_i$, then $\vF_{\vT_1} = \vF_{\vT_2}$. 
\item 
Let $\theta\in \scr C(\frak a,\beta)$ have endo-class $\vT$. The character $\check\theta = \theta^{-1}$ of $H^1(\beta,\frak a)$ lies in $\scr C(\frak a,-\beta)$ and its endo-class $\vT\spcheck$ satisfies $\vF_{\vT\spcheck} = \vF_\vT$. 
\endroster 
\endproclaim 
\demo{Proof} Both assertions follow directly from the Observation concluding 4.3. \qed 
\enddemo 
\subhead 
4.5 
\endsubhead 
The functions $\vF_\vT$ reflect the approximation properties intrinsic to the concept of endo-class. 
\proclaim{Proposition} 
For $i=1,2$, let $[\frak a,m_i,0,\beta_i]$ be a simple stratum in $\End FV$, for a finite-dimensional $F$-vector space $V$. Let $\theta_i\in \scr C(\frak a,\beta_i)$ and let $\vT_i$ be the endo-class of $\theta_i$. If $t\ge 0$ is an integer such that the restrictions $\theta_i\,\big|\,{H^{1+t}(\beta_i,\frak a)}$ intertwine in $\Aut FV$, then 
$$
\vF_{\vT_1}(x) = \vF_{\vT_2}(x), \quad x\ge t/e_\frak a. 
$$ 
\endproclaim 
\demo{Proof} 
Choose a simple stratum $[\frak a,m_i,t,\gamma_i]$ equivalent to $[\frak a,m_i,t,\beta_i]$. In particular, $H^{1+t}(\beta_i,\frak a) = H^{1+t}(\gamma_i,\frak a)$ and the character $\theta_i^t = \theta_i\,\big|\,{H^{1+t}(\beta_i,\frak a)}$ lies in $\scr C(\frak a,t,\gamma_i)$. The truncated simple characters $\theta_i^t$ intertwine and and so are conjugate in $\Aut FV$ \cite{9} (3.5.11). By \cite{9} (3.5.1), $k_0(\gamma_1,\frak a) = k_0(\gamma_2,\frak a)$. The proposition now follows from the definition 4.4.1, the Observation of 4.3 and induction along the stratum $[\frak a,m_1,0,\beta_1]$. \qed 
\enddemo 
\subhead 
4.6 
\endsubhead 
Let $\vT \in \Scr E(F)$, and let $K/F$ be a finite, tamely ramified field extension. Let $\vT_i^K \in \Scr E(K)$, $1\le i\le r$, be the set of {\it $K/F$-lifts\/} of $\vT$ \cite{2}. If $\theta\in \scr C(\frak a,\beta)$ is a realization of $\vT$, the $\vT^K_i$ are in bijection with the simple components of the semisimple $K$-algebra $K\otimes_F F[\beta]$. The relation between $\vF_\vT$ and the functions $\vF_{\vT^K_i}$ is, in general, quite intricate but we shall only need a special case. 
\par 
Say that $\vT \in \Scr E(F)$ is {\it totally wild\/} if $e_\vT = \deg\vT = p^r$, for an integer $r\ge0$. If $\theta\in \scr C(\frak a,\beta)$ is a realization of $\vT$, then $\vT$ is totally wild if and only if the field extension $F[\beta]/F$ is totally wildly ramified. 
\proclaim{Proposition} 
Let $\vT \in \Scr E(F)$ be totally wild. If $K/F$ is a finite, tamely ramified field extension, then $\vT$ has a unique $K/F$-lift $\vT^K$. If $e= e(K|F)$, then 
$$
\vF_{\vT^K}(x) = e\,\vF_\vT(x/e), \quad x\ge0. 
$$ 
\endproclaim 
\demo{Proof} 
Let $\deg\vT = p^a$, $a\ge1$. If $V$ is an $F$-vector space of dimension $p^r$, there is a simple stratum $[\frak a_0,m,0,\beta]$ in $\End FV$ such that $\scr C(\frak a_0,\beta)$ contains a character $\theta_0$ of endo-class $\vT$. If $E = F[\beta]$ then $E/F$ is totally ramified of degree $p^r$. Having identified $E$, we may as well take $V = E$, and then $\frak a_0$ is the unique $E$-pure hereditary $\frak o_F$-order in $A(E) = \End FE$. 
\par 
The algebra $K\otimes_FE$ is a field, which we denote $KE$. In particular, $\vT$ admits a unique $K/F$-lift $\vT^K$. 
\par 
Let $A = \End F{KE}$ and let $\frak a$ be the unique $KE$-pure hereditary $\frak o_F$-order in $A$. The quadruple $[\frak a,em,0,\beta]$ is a simple stratum in $A$ and there is a simple character $\theta\in \scr C(\frak a,\beta)$ of endo-class $\vT$. As $e_\frak a = ep^r$, so 
$$ 
\scr S_\vT = \{x/ep^r:x\in \roman S_{[\frak a,\beta]}\}, 
$$ 
in the notation of 4.3. 
\par 
Let $B = \End K{KE}$ be the $A$-centralizer of $K$ and $\frak b  = \frak a\cap B$. Thus $\frak b$ is the unique $KE$-pure hereditary $\frak o_K$-order in $B$. The stratum $[\frak b,em,0,\beta]$ is simple \cite{2} (2.4). Further, $H^1(\beta,\frak b) = H^1(\beta,\frak a)\cap B^\times$, and the character $\theta^K = \theta\,\big|\,{H^1(\beta,\frak b)}$ lies in $\scr C(\frak b,\beta,\psi_K)$, where $\psi_K = \psi_F\circ \roman{Tr}_{K/F}$ \cite{2} (7.7). The endo-class of $\theta^K$ over $K$ is then $\vT^K$. 
\proclaim{Lemma} 
The sets $\roman S_{[\frak b,\beta]}$, $\roman S_{[\frak a,\beta]}$ are equal. 
\endproclaim 
\demo{Proof} 
We proceed by induction along $\beta$, in the manner of many proofs in \cite{9}. Suppose first that $\beta$ is minimal over $F$. It is then minimal over $K$ \cite{2} (2.4) and the field extensions $F[\beta]/F$, $K[\beta]/K$ are totally ramified of the same degree. The lemma holds in this case. 
\par
We therefore assume $r = -k_0(\beta,\frak a) < em$. Let $s = -k_0(\beta,\frak b)$. According to \cite{2} (2.4), we have $s\ge r$. We show that $s = r$ in this case: assume for a contradiction that $s>r$. We choose a simple stratum $[\frak b,em,s{-}1,\gamma]$ in $B$, equivalent to $[\frak b,em,s{-}1,\beta]$, such that $[\frak a,em,s{-}1,\gamma]$ is simple: this we may do by \cite{2} (3.8). Certainly $[\frak a,em,s{-}1,\gamma]$ is equivalent to $[\frak a,em,s{-}1,\beta]$, which is not simple. It follows from \cite{9} (2.4.1) that $F[\gamma]/F$ is totally wildly ramified and $[F[\gamma]{:}F] < [F[\beta]{:}F]$. Thus $[K[\gamma]{:}K] < [K[\beta]{:}K]$, implying that $[\frak b,em,s{-}1,\beta]$ is not simple. This eliminates the possibility $s>r$. 
\par
We conclude that the sets $\roman S_{[\frak a,\beta]}$, $\roman S_{[\frak b,\beta]}$ have the same least element $r =-k_0(\beta,\frak a)$. By definition, $\roman S_{[\frak a,\beta]} = \{r\}\cup \roman S_{[\frak a,\gamma]}$ and $\roman S_{[\frak b,\beta]} = \{r\}\cup \roman S_{[\frak b,\gamma]}$. Inductively, $\roman S_{[\frak a,\gamma]} = \roman S_{[\frak b,\gamma]}$, and the lemma is proved. \qed 
\enddemo 
We deduce that 
$$ 
\scr S_{\vT^K} = \{ey: y\in \scr S_\vT\}. 
\tag 4.6.1 
$$ 
The proof of the lemma also shows that, if $x>0$ and $x\notin \scr S_\vT$, then 
$$ 
(KE)_{ex} = KE_x \quad \text{whence} \quad d_\vT(x) = d_{\vT^K}(ex),  
\tag 4.6.2 
$$ 
in the notation of 4.3. 
\par 
Set $\phi(x) = e\vF_\vT(x/e)$. The functions $\phi$ and $\vF_{\vT^K}$ are continuous, and smooth outside of $\scr S_{\vT^K}$. Also, by (4.6.2), $\phi'(x) = \vF_\vT'(x)$ for $x\notin \scr S_{\vT^K}$. In other words, $\phi(x) = \vF_{\vT^K}(x){+}c$, for a constant $c$. However, for $x$ sufficiently large, $\phi(x) = \vF_{\vT^K}(x) = x$, so $c=0$ as required to prove the proposition. \qed 
\enddemo 
\head\Rm 
5. Ultrametric on simple characters 
\endhead 
We re-examine the conductor formula of \cite{8}, interpreting it in terms of the structure functions $\vF_\vT$ of 4.4 and a canonical ultrametric on the set $\Scr E(F)$ of endo-classes of simple characters over $F$. 
\subhead 
5.1 
\endsubhead 
Let $\vT_1$, $\vT_2 \in \Scr E(F)$, $\vT_i \neq \bk0$. There is an $F$-vector space $V$ of finite dimension, and a hereditary order $\frak a$ in $\End FV$, such that $\frak a$ carries realizations of both $\vT_i$. That is, there are simple strata $[\frak a,m_i,0,\beta_i]$ in $\End FV$ and simple characters $\theta_i\in \scr C(\frak a,\beta_i)$ such that $\theta_i$ is of endo-class $\vT_i$. 
\par 
Let $l\ge0$ be the least integer such that the characters $\theta_i \,\big|\,{H^{l+1}(\beta_i,\frak a)}$ intertwine (and are therefore conjugate \cite{9} (3.5.11)) in $\Aut FV$. We define 
$$
\Bbb A(\vT_1,\vT_2) = \Bbb A(\vT_2,\vT_1) = l/e_\frak a. 
\tag 5.1.1 
$$ 
The definition is independent of all choices: see the discussion in \cite{8} 6.15. One may treat the trivial class $\bk0$ on the same basis, but it is quicker to simply define 
$$
\Bbb A(\vT,\bk 0) = m_\vT,\quad \vT \in \Scr E(F). 
\tag 5.1.2 
$$ 
\proclaim{Proposition} 
\roster 
\item 
Let $\vT,\vT'\in \Scr E(F)$. If $m_\vT < m_{\vT'}$, then $\Bbb A(\vT,\vT') = m_{\vT'}$. 
\item 
If $\vT_1,\vT_2\in \Scr E(F)$, then $\Bbb A(\vT_1,\vT_2) = 0$ if and only if $\vT_1 = \vT_2$. 
\item 
If $\vT_1,\vT_2,\vT_3 \in \Scr E(F)$, then 
$$
\Bbb A(\vT_1,\vT_3) \le \roman{max}\,\{\Bbb A(\vT_1,\vT_2), \Bbb A(\vT_2,\vT_3)\}. 
\tag 5.1.3 
$$ 
\endroster 
\endproclaim 
\demo{Proof} 
Part (1) follows from \cite{9} (2.6.3). In (2), we find a hereditary order $\frak a$ in some $A = \End FV$, a simple stratum $[\frak a,m_i,0,\beta_i]$ and a simple character $\theta_i\in \scr C(\frak a,0,\beta_i)$ of endo-class $\vT_i$, $i=1,2$. The assertion $\Bbb A(\vT_1,\vT_2) = 0$ is equivalent to the characters $\theta_i$ of $H^1(\beta_i,\frak a)$ being conjugate in $\Aut FV$. This, in turn, is equivalent to $\vT_1 = \vT_2$. 
\par 
In (3), we may take simultaneous realizations $\theta_i\in \scr C(\frak a,\beta_i)$ of $\vT_i$, $i=1,2,3$, in some $G = \Aut FV$. Let $t_{ij}$ be the least non-negative integer such that $\theta_i\,\big|\,{H^{1+t_{ij}}(\beta_i,\frak a)}$ is $G$-conjugate to $\theta_j\,\big|\,{H^{1+t_{ij}}(\beta_i,\frak a)}$. Thus $\Bbb A(\vT_i,\vT_j) = t_{ij}/e_\frak a$. By symmetry, we may assume that $t_{12}\le t_{23}$. Replacing the $\theta_i$ by conjugates, we may further assume that 
$$ 
\align 
H^{1+t_{12}}(\beta_1,\frak a) &= H^{1+t_{12}}(\beta_2,\frak a), \\ 
H^{1+t_{23}}(\beta_2,\frak a) &= H^{1+t_{23}}(\beta_3,\frak a), 
\endalign 
$$ 
and that 
$$ 
\alignedat3 
\theta_1(g) &= \theta_2(g), &\quad &g\in H^{1+t_{12}}(\beta_1,\frak a), \\ 
\theta_2(h) &= \theta_3(h), &\quad &h\in H^{1+t_{23}}(\beta_2,\frak a) = H^{1+t_{23}}(\beta_1,\frak a). 
\endalignedat 
$$ 
Thus $\theta_1$ agrees with $\theta_3$ on $H^{1+t_{23}}(\beta_i,\frak a)$. It follows that $t_{13} \le t_{23}$, as required to prove (5.1.3). \qed 
\enddemo 
In summary, the pairing $\Bbb A$ defines an {\it ultrametric\/} on the set $\Scr E(F)$. It is natural to re-state 4.5 Proposition in terms of $\Bbb A$. 
\proclaim{Corollary} 
If $\vT, \vT' \in \Scr E(F)$, then $\vF_\vT(a) = \vF_{\vT'}(a)$ for all $a\ge \Bbb A(\vT,\vT')$. 
\endproclaim 
\subhead 
5.2 
\endsubhead 
In some circumstances, a different language is clearer. Let $[\frak a,m,0, \allowmathbreak \beta]$ be a simple stratum in some matrix algebra $A = \End FV$, and let $\theta\in\scr C(\frak a,\beta)$. Let $\eps > 0$ and let $t$ be the greatest integer such that $t/e_\frak a<\eps$. In particular, $t\ge 0$. The {\it $\eps$-truncation of\/} $\theta$, denoted  $\roman{tc}_\eps(\theta)$, is the character $\theta\,\big|\,H^{1+t}(\beta,\frak a)$. 
\par 
Using the general machinery of \cite{2} section 8, we may form the endo-class of $\roman{tc}_\eps(\theta)$: if $\vT$ is the endo-class of $\theta$, we denote the endo-class of $\roman{tc}_\eps(\theta)$ by $\roman{tc}_\eps(\vT)$. This depends only on $\vT$ and $\eps$. The definition of the ultrametric $\Bbb A$ then implies: 
\proclaim{Proposition} 
Let $\eps > 0$. If $\vT_1,\vT_2 \in \Scr E(F)$, then $\Bbb A(\vT_1,\vT_2) < \eps$ if and only if $\roman{tc}_\eps(\vT_1) = \roman{tc}_\eps (\vT_2)$. 
\endproclaim 
\subhead 
5.3 
\endsubhead 
The following, more delicate, property is needed in certain situations. 
\proclaim{Density Lemma} 
Let $\vT \in \Scr E(F)$. The set $\{\Bbb A(\vT,\vX): \vX \in \Scr E(F)\}$ is dense in the half line $\{x\ge0: x\in \Bbb R\}$. 
\endproclaim 
\demo{Proof} 
Let $x\in \Bbb Q$, $x>0$. If $x> m_\vT$, there exists $\vX\in \Scr E(F)$ such that $m_\vX = x$, by 4.2 Proposition. This gives $\Bbb A(\vT,\vX) = x$, by 5.1 Proposition, so it is enough to treat the case $x < m_\vT$. 
\proclaim{Lemma} 
Let $[\frak a,m,0,\beta]$ be a simple stratum in a matrix algebra $A = \M nF$ and let $\theta\in \scr C(\frak a,\beta)$. Let $k$  be an integer, $1\le k \le m$. There exists a simple stratum $[\frak a,m,0,\beta']$ in $A$ and $\theta'\in \scr C(\frak a,\beta')$ such that 
\roster 
\item $H^k(\beta',\frak a) = H^k(\beta,\frak a)$, 
\item $\theta'$ agrees with $\theta$ on $H^{k+1}(\beta,\frak a)$, and 
\item the characters $\theta$, $\theta'$ do not intertwine on $H^k(\beta,\frak a)$. 
\endroster 
\endproclaim 
\demo{Proof} 
We first reduce to the case in which the stratum $[\frak a,m,k{-}1,\beta]$ is simple. Suppose it is not. We choose a simple stratum $[\frak a,m,k{-}1,\gamma]$ equivalent to $[\frak a,m,k{-}1,\beta]$. Directly from the definitions in \cite{9} Chapter 3 we have 
$$
H^k(\gamma,\frak a) = H^k(\beta,\frak a), \quad \scr C(\frak a,k{-}1,\gamma) = \scr C(\frak a, k{-}1, \beta). 
$$ 
In particular, there exists $\xi\in \scr C(\frak a,\gamma)$ agreeing with $\theta$ on $H^k(\beta,\frak a)$. We may now work with the pair $(\gamma,\xi)$ in place of $(\beta,\theta)$. 
\par 
We revert to our original notation, assuming that $[\frak a,m,k{-}1,\beta]$ is simple. Let $B$ denote the $A$-centralizer of $\beta$ and let $\frak b = \frak a\cap B$. We choose a simple stratum $[\frak b,k,k{-}1,\alpha]$ in $B$. Writing $\frak p = \roman{rad}\,\frak a$, let $a\in \frak p^{-k}$ satisfy $s_\beta(a) = \alpha$, where $s_\beta:A \to B$ is a tame corestriction relative to $F[\beta]/F$. The stratum $[\frak a,m,k{-}1, \allowmathbreak \beta{+}a]$ is then equivalent to a simple stratum $[\frak a,m,k{-}1,\beta']$ \cite{9} (2.2.3). Let $\psi_a$ denote the character 
$$
1{+}x \longmapsto \psi_F(\roman{tr}_A(ax)), \quad x\in \frak p^k, 
$$ 
of $U^k_\frak a$. The character $\theta' = \theta\psi_a$ of $H^k(\beta',\frak a) = H^k(\beta,\frak a)$ then lies in $\scr C(\frak a,k{-}1,\beta')$ and agrees with $\theta$ on $H^{1+k}(\beta,\frak a)$. However, 2.8 Proposition of \cite{11} implies that the characters $\theta$, $\theta'$ of $H^k(\beta,\frak a)$ do not intertwine. \qed 
\enddemo 
In the context of the lemma, let $\theta$, $\theta'$ have endo-class $\vT$, $\vT'$ respectively. Thus $\Bbb A(\vT,\vT') = k/e_\frak a$. The only restrictions on the rational number $k/e_\frak a$ are that $e_\frak a$ be divisible by $e_\vT$ and $k/e_\frak a \le m_\vT$. Such values are dense in the region $0<x<m_\vT$. \qed 
\enddemo 
\subhead 
5.4 
\endsubhead 
We recall the notation of the introduction: $\Ao nF$ is the set of equivalence classes of irreducible cuspidal representations of $\GL nF$. We set $\wG F = \bigcup_{n\ge 1} \Ao nF$ and, for $\pi \in \wG F$, we write $\roman{gr}(\pi) = n$ to indicate $\pi \in \Ao nF$. Let $\check\pi$ be the contragredient of $\pi$. 
\par 
A representation $\pi \in \wG F$ contains a simple character $\theta_\pi$. The conjugacy class of $\theta_\pi$ in $\GL nF$ is determined uniquely by $\pi$: see, for instance, Corollary 1 of \cite{6}. In particular, $\pi$ determines the endo-class $\vT = \vT_\pi$ of $\theta$. 
\par 
We recall the definition of the {\it Swan exponent\/} $\sw(\pi_1\times\pi_2)$ of a pair of representations $\pi_1,\pi_2 \in \wG F$. Set $n_i = \roman{gr}(\pi_i)$, let $\psi$ be a non-trivial character of $F$, let $s$ be a complex variable and $q$ the cardinality of the residue class field of $F$. Let $\ve(\pi_1\times\pi_2,s,\psi)$ be the Rankin-Selberg local constant of \cite{18} and \cite{23}. This is a monomial in $q^{-s}$ of degree $n_1n_2c(\psi){+}\roman{Ar}(\pi_1\times\pi_2)$, where $c(\psi)$ is an integer depending only on $\psi$, and the {\it Rankin-Selberg exponent\/} $\roman{Ar}(\pi_1\times\pi_2)$ is an integer depending only on the $\pi_i$. Define an integer $d(\pi_1,\pi_2)$ as the number of unramified characters $\chi$ of $F^\times$ such that $\chi\pi_1 \cong \check\pi_2$. In particular, $d(\pi_1,\pi_2) = 0$ if $n_1\neq n_2$. The Swan exponent is then 
$$ 
\sw(\pi_1\times\pi_2) = \roman{Ar}(\pi_1\times\pi_2) - n_1n_2 + d(\pi_1,\pi_2). 
$$ 
Reformulating 6.5 Theorem of \cite{8} in our present notation, we find: 
\proclaim{Conductor formula} 
For $i=1,2$, let $\pi_i \in \wG F$ and set $\vT_i = \vT_{\pi_i}$. If $a = \Bbb A(\vT_1,\vT_2)$, then 
$$
\frac{\roman{sw}(\check\pi_1\times \pi_2)}{\roman{gr}(\pi_1)\,\roman{gr}(\pi_2)} = \vF_{\vT_1}(a) = \vF_{\vT_2}(a). 
\tag 5.4.1 
$$ 
\endproclaim 
If we take $\pi\in \Ao nF$ and let $\iota$ be the trivial character of $\GL1F$, we get the special case ({\it cf\.} (5.1.2)) 
$$ 
\sw(\pi\times\iota)/n = \sw(\pi)/n = m_{\vT_\pi}. 
\tag 5.4.2 
$$ 
\proclaim{Proposition} 
Let $\vT \in \Scr E(F)$, and let $\pi \in \wG F$ satisfy $\vT_\pi = \vT$. The function $\vF_\vT$ is the unique continuous function on the positive real axis such that 
$$
\frac{\roman{sw}(\check\pi\times \rho)}{\roman{gr}(\pi)\,\roman{gr}(\rho)} = \vF_\vT\big(\Bbb A(\vT, \vT_\rho)\big), 
$$ 
for all $\rho \in \wG F$. 
\endproclaim 
\demo{Proof} 
This follows from (5.4.1), the continuity of the function $\vF_\vT$ (4.4 Proposition) and the Density Lemma of 5.3. \qed 
\enddemo 
The proposition has a consequence which is useful when making more general conductor estimates, although we do not need it here. For $i=1,2$, let $\vT_i\in \Scr E(F)$ and choose $\pi_i\in \wG F$ such that $\vT_i = \vT_{\pi_i}$. Let $\roman{gr}(\pi_i) = n_i$. The quantity 
$$ 
\vs(\pi_1,\pi_2) = \sw(\check\pi_1\times\pi_2)/n_1n_2 
$$ 
depends only on the $\vT_i$, not on the choices of $\pi_i$: this is a consequence of the proposition. We therefore write $\vs(\vT_1,\vT_2) = \vs(\pi_1,\pi_2)$. 
\proclaim{Corollary} 
The pairing $\vs$ on the set $\Scr E(F)$ satisfies the ultrametric inequality: if $\vT_1,\vT_2,\vT_3 \in \Scr E(F)$, then 
$$
\vs(\vT_1,\vT_2) \le \roman{max}\,\{\vs(\vT_1,\vT_3),\vs(\vT_3,\vT_2)\}. 
$$ 
\endproclaim 
\demo{Proof} 
The pairing $\vs$ is symmetric: $\vs(\vT_1,\vT_2) = \vs(\vT_2,\vT_1)$. We may assume, by symmetry, that $\Bbb A(\vT_1,\vT_3) \le \Bbb A(\vT_3,\vT_2)$. The function $\vF_{\vT_3}$ is increasing, so 
$$ 
\vs(\vT_1,\vT_3) = \vF_{\vT_3}(\Bbb A(\vT_1,\vT_3)) \le \vF_{\vT_3}(\vT_2,\vT_3)) = \vs(\vT_2,\vT_3). 
$$ 
We are thus reduced to checking that $\vs(\vT_1,\vT_2) \le \vs(\vT_2,\vT_3)$. However, the ultrametric inequality for $\Bbb A$ and our hypothesis give $\Bbb A(\vT_1,\vT_2) \le \Bbb A(\vT_3,\vT_2)$ so
$$ 
\vs(\vT_1,\vT_2) = \vF_{\vT_2}(\vT_1,\vT_2) \le \vF_{\vT_2}(\vT_3,\vT_2) = \vs(\vT_3,\vT_2), 
$$ 
as required. \qed 
\enddemo 
\subhead 
5.5 
\endsubhead 
We give a property of the ultrametric $\Bbb A$ relative to tame lifting, as in \cite{2} (see also 4.6 above). For clarity, we temporarily write $\Bbb A_F$ for the canonical ultrametric on $\Scr E(F)$ and $\Bbb A_K$ for that on $\Scr E(K)$, where $K/F$ is a finite tame extension. 
\proclaim{Proposition} 
Let $\vT, \vU \in \Scr E(F)$ and let $K/F$ be a finite tame extension with $e(K|F) = e$. If $\vT_i$, $1\le i\le r$ are the $K/F$-lifts of $\vT$ and $\vU_j$, $1\le j\le s$, those of $\vU$, then 
$$ 
e\,\Bbb A_F(\vT,\vU) =  \underset{i,j} \to{\roman{min}}\,\Bbb A_K(\vT_i,\vU_j) = \underset{j} \to{\roman{min}}\,\Bbb A_K(\vT_1,\vU_j) . 
$$ 
\endproclaim 
\demo{Proof} 
From (9.8) Theorem of \cite{2} we deduce that $e\,\Bbb A_F(\vT,\vU) \le \Bbb A_K(\vT_i,\vU_j)$, for all $i$ and $j$. On the other hand, \cite{2} (9.12) Corollary implies that, for any $i$, there exists $j$ such that $e\,\Bbb A_F(\vT,\vU) \ge \Bbb A_K(\vT_i,\vU_j)$, whence the result follows. \qed 
\enddemo 
\head \Rm 
6. Comparison via the Langlands correspondence 
\endhead 
We use the local Langlands correspondence to connect the preceding lines of thought. 
\subhead 
6.1 
\endsubhead 
We recall formally some matters mentioned in introduction. Using the notation of 5.4, the Langlands correspondence is a canonical bijection 
$$ 
\align 
\wG F &\longrightarrow \wW F, \\ \pi &\longmapsto \upr L\pi, 
\endalign 
$$ 
with, among others, the following properties: 
$$ 
\aligned 
\dim\upr L\pi &= \roman{gr}(\pi), \\ \upr L(\check\pi) &= \big(\upr L\pi\big)\spcheck, \\ 
\ve(\pi\times\rho,s,\psi) &= \ve(\upr L\pi\otimes \upr L\rho,s,\psi), 
\endaligned 
\qquad \pi,\rho \in \wG F. 
\tag 6.1.1 
$$ 
Here, the second $\ve$ is the Langlands-Deligne local constant. The correspondence also respects twisting with characters. The definition of $\sw(\pi_1\times\pi_2)$ in 5.4 thus implies  
$$ 
\sw(\pi\times\rho) = \sw(\upr L\pi\otimes\upr L\rho), \quad \pi,\rho \in \wG F. 
$$ 
We prefer to write $\vs(\pi_1\times\pi_2) = \sw(\pi_1\times\pi_2)/\roman{gr}(\pi_1)\roman{gr}(\pi_2)$, so that 
$$ 
\vs(\check\pi\times \rho) = \vs(\upr L{\check\pi} \otimes \upr L\rho), \quad \pi,\rho \in \wG F. 
\tag 6.1.2 
$$ 
A representation $\pi \in \wG F$ determines an endo-class $\vT_\pi \in \Scr E(F)$, as recalled in 5.4. On the other hand, a representation $\sigma\in \wW F$ determines an orbit $[\sigma;0]^+ \in \scr W_F\backslash \wP F$, as in 2.2. 
\proclaim{First ramification theorem} 
Let $\vT \in \Scr E(F)$ and choose $\pi \in \wG F$ such that $\vT_\pi = \vT$. The conjugacy class $\upr L\vT = [\upr L\pi;0]^+ \in \scr W_F\backslash \wP F$ depends only on $\vT$ and not on the choice of $\pi$. The map 
$$ 
\aligned 
\Scr E(F) &\longrightarrow \scr W_F\backslash \wP F, \\ 
\vT &\longmapsto \upr L\vT, 
\endaligned 
$$ 
is a canonical bijection. 
\endproclaim 
\demo{Proof} 
See \cite{4} 8.2 Theorem, \cite{7} 6.1. \qed 
\enddemo  
\subhead 
6.2 
\endsubhead 
Let $\vT \in \Scr E(F)$. Choose $\pi \in \wG F$ so that $\vT_\pi = \vT$ and write $\upr L\pi = \sigma$. The decomposition function $\vS_\sigma$ depends only on $[\sigma;0]^+ = \upr L\vT$, so we use the notation $\vS_\sigma = \vS_{\upr L\vT}$. Combining (6.1.2) with (3.1.4) and 5.4 Proposition, we find 
$$ 
\aligned 
\vF_\vT(\Bbb A(\vT, \vT_\rho)) &= \vs(\check\pi\times \rho) = \vs(\upr L{\check\pi} \otimes \upr L\rho)\\ 
&= \vS_{\upr L\vT}(\Delta(\upr L\vT, \upr L\vT_\rho)), 
\endaligned \quad \rho\in \wG F. 
$$ 
In other words, 
$$ 
\vF_\vT(\Bbb A(\vT,\vU)) = \vS_{\upr L\vT}(\Delta(\upr L\vT,\upr L\vU)), \quad \vT,\vU \in \Scr E(F). 
\tag 6.2.1 
$$ 
We accordingly define the {\it Herbrand function\/} $\Psi_\vT$ of $\vT$ by 
$$ 
\Psi_\vT = \vF_\vT^{-1}\circ \vS_{\upr L\vT}, \quad \vT \in \Scr E(F). 
\tag 6.2.2 
$$ 
\proclaim{Proposition} 
Let $\vT \in \Scr E(F)$. 
\roster 
\item 
The function $\Psi_\vT$ is continuous, strictly increasing and piece-wise linear in the region $x\ge0$. It is smooth except at a finite set of points.  
\item 
It satisfies $\Psi_\vT(0) = 0$ and $\Psi_\vT(x) = x$ for $x\ge m_\vT$. 
\endroster 
\endproclaim 
\demo{Proof} 
Part (1) combines 4.4 Proposition with 3.1 Proposition. In part (2), we choose $\pi \in \wG F$ such that $\vT_\pi = \vT$ and set $\sigma = \upr L\pi$. Thus 
$$
\vF_\vT(0) = \vs(\check\pi\times \pi) = \vs(\check\sigma \otimes \sigma) = \vS_{\upr L\vT}(0), 
$$ 
whence $\Psi_\vT(0) = 0$. By (5.4.2), $m_\vT = \vs(\pi) = \vs(\sigma) = \roman{sl}(\sigma)$, so the second assertion in (2) follows from (3.1.2) and (4.4.1). \qed 
\enddemo 
\subhead 
6.3 
\endsubhead 
We pause, to tie up some loose ends. Since $\sw(\pi) = \sw(\upr L\pi)$ and $\roman{gr}(\pi) = \dim\upr L\pi$,  2.4 Proposition follows from 4.2 Proposition. The Density Lemma of 3.3 follows from that of 5.3 and the continuity of the strictly increasing function $\Psi_\vT^{-1}$. This proves 3.3 Proposition.  
\subhead 
6.4 
\endsubhead 
We prove our first main result. 
\proclaim{Higher Ramification Theorem} 
Let $\vT \in \Scr E(F)$, let $\eps>0$ and $\delta = \Psi_\vT(\eps)$. If $\vU \in \Scr E(F)$ then 
$$ 
\alignedat3 
\Delta(\upr L\vT,\upr L\vU) &< \eps \quad &\Longleftrightarrow \quad \Bbb A(\vT,\vU) &< \delta, \\ 
\Delta(\upr L\vT,\upr L\vU) &\le \eps \quad &\Longleftrightarrow \quad \Bbb A(\vT,\vU) &\le \delta. 
\endalignedat 
$$
\endproclaim 
\demo{Proof} 
Let $\vT \in \Scr E(F)$ and $\delta > 0$. The endo-class $\vT$ determines the function $\vF_\vT$ and the orbit $\upr L\vT \in \scr W_F\backslash \wP F$, whence it determines the function $\vS_{\upr L\vT}$. For $\vU \in \Scr E(F)$, the condition $\Bbb A(\vT,\vU) < \delta$ implies 
$$ 
\Delta(\upr L\vT, \upr L\vU) = \vS_{\upr L\vT}^{-1}\vF_\vT(\Bbb A(\vT,\vU)) = \Psi_\vT^{-1}(\Bbb A(\vT,\vU)) < \Psi^{-1}_\vT(\delta), 
$$ 
since the function $\Psi_\vT$ is strictly increasing (6.2 Proposition). Indeed, the converse holds for the same reason: if $\Delta(\upr L\vT,\upr L\vU) < \Psi^{-1}_\vT(\delta)$, then $\Bbb A(\vT,\vU) < \delta$. The same argument proves the second assertion. \qed 
\enddemo 
\subhead 
6.5 
\endsubhead 
We give a more concrete variant of the main theorem. We first need a technical result.
\proclaim{Lemma} 
If $\vT,\vU \in \Scr E(F)$ and $x\ge \Bbb A(\vT,\vU)$, then $\Psi^{-1}_\vT(x) = \Psi^{-1}_\vU(x)$. 
\endproclaim 
\demo{Proof} 
Let $\delta > \Bbb A(\vT,\vU)$ and set $\eps = \Psi_\vT^{-1}(\delta)$. An endo-class $\vX\in \Scr E(F)$ thus satisfies $\Bbb A(\vX,\vU) < \delta$ if and only if $\Bbb A(\vX,\vT) <\delta$. The second condition is equivalent to $\Delta(\upr L\vX,\upr L\vT) < \eps$ by the theorem, while the first is equivalent to $\Delta(\upr L\vU,\upr L\vX) < \Psi^{-1}_\vU(\delta)$. On the other hand, 
$$ 
\Delta(\upr L\vU,\upr L\vX) \le \roman{max}\big\{\Delta(\upr L\vU,\upr L\vT), \Delta(\upr L\vT,\upr L\vX)\big\} 
< \eps. 
$$
It follows that $\Psi_\vU^{-1}(\delta) \le \eps = \Psi_\vT^{-1}(\delta)$ for $\delta > \Bbb A(\vT,\vU)$. By symmetry, 
$$ 
\Psi_\vU^{-1}(\delta)  = \Psi_\vT^{-1}(\delta), \quad \delta > \Bbb A(\vT,\vU). 
$$ 
By continuity, the relation holds for $ \delta \ge \Bbb A(\vT,\vU)$. \qed 
\enddemo 
\remark{Remark} 
Under the hypotheses of the lemma, we may equally deduce that $\Psi_\vT(y) = \Psi_\vU(y)$ when $y\ge \Delta(\upr L\vT,\upr L \vU)$. 
\endremark 
We now use the notation of 5.2 for truncated endo-classes. 
\proclaim{Corollary} 
\roster 
\item 
Let $\vT \in \Scr E(F)$ and $\delta>0$. There is a unique pair $(\eps,\xi)$, where $\eps>0$ and $\xi\in \scr W_F\backslash \wR F\eps$, with the following property: a representation $\pi \in \wG F$ satisfies $\roman{tc}_\delta(\vT_\pi) = \roman{tc}_\delta(\vT)$ if and only if the representation $\xi$ is equivalent to a component of $\upr L\pi \,\big|\,\scr R_F(\eps)$. 
\item 
Let $\eps>0$ and $\xi\in \scr W_F\backslash \wR F\eps$. There exist $\vT \in \Scr E(F)$ and $\delta>0$ with the following property: a representation $\pi \in \wG F$ satisfies $\roman{tc}_\delta(\vT_\pi) = \roman{tc}_\delta(\vT)$ if and only if the representation $\xi$ is equivalent to a component of $\upr L\pi \,\big|\,\scr R_F(\eps)$. The pair $(\eps,\xi)$ determines the truncated endo-class $\roman{tc}_\delta(\vT)$ uniquely. 
\endroster 
\endproclaim 
\demo{Proof} 
It is enough to treat part (1). Set $\eps = \Psi_\vT^{-1}(\delta)$ and let $\xi$ be the conjugacy class of an irreducible component of $\upr L\vT$ on $\scr R_F(\eps)$. From the theorem, a class $\vU\in \Scr E(F)$ satisfies $\Bbb A(\vU,\vT) <\delta$ if and only if $\Delta(\upr L\vU,\upr L\vT) <\eps$. The first of these conditions is equivalent to $\roman{tc}_\delta(\vU) = \roman{tc}_\delta(\vT)$ (5.2 Proposition) while the second is equivalent to $\upr L\vU$ containing $\xi$, by the definition of $\Delta$. All assertions now follow. \qed 
\enddemo 
\head\Rm 
7. The Herbrand function of an endo-class 
\endhead 
We give a procedure for determining the Herbrand function $\Psi_\vT$ of an endo-class $\vT \in \Scr E(F)$. 
\subhead 
7.1 
\endsubhead 
Fundamental to the method is the following lifting property. 
\proclaim{Proposition} 
Let $K/F$ be a finite, tame extension and set $e(K|F) = e$. If $\vT\in \Scr E(F)$ and if $\vT^K\in \Scr E(K)$ is a $K/F$-lift of $\vT$, then 
$$ 
\Psi_{\vT^K}(x) = e\,\Psi_\vT(e^{-1}x), \quad x\ge 0. 
\tag 7.1.1 
$$ 
\endproclaim 
\demo{Proof} 
Using transitivity of tame lifting, we reduce immediately to the case where the tame extension $K/F$ is Galois. Write $\vG = \Gal KF$ and let $\vU\in \Scr E(F)$. Let $\vU^K$ be a $K/F$-lift of $\vU$. Write $\Bbb A_F$, $\Bbb A_K$ for the canonical ultrametrics on $\Scr E(F)$, $\Scr E(K)$ respectively. We choose the lift $\vU^K$ so that 
$$ 
\Bbb A_K(\vT^K,\vU^K) \le \Bbb A_K(\vT^K,\gamma\vU^K), \quad \gamma \in \vG. 
$$ 
The function $\Psi_{\vT^K}$ is strictly increasing, so writing $\Delta_F$, $\Delta_K$ for the canonical ultrametrics on $\scr W_F\backslash \wP F$, 
$\scr W_K\backslash \wP K$ respectively, we have 
$$ 
\Delta_K(\upr L\vT^K,\upr L\vU^K) \le \Delta_K(\upr L\vT^K, \upr L(\gamma\vU^K)), \quad \gamma\in \vG. 
$$
The canonical bijection $\Scr E(K) \to \scr W_K\backslash \wP K$ is $\vG$-equivariant, so this reads 
$$ 
\Delta_K(\upr L\vT^K,\upr L\vU^K) \le \Delta_K(\upr L\vT^K, \gamma\,\upr L\vU^K), \quad \gamma\in \vG, 
$$ 
whence 
$$ 
\align 
\Bbb A_F(\vT,\vU) &= e^{-1}\Bbb A_K(\vT^K,\vU^K), \\ 
\Delta_F(\upr L\vT,\upr L\vU) &= e^{-1}\Delta_K(\upr L\vT^K,\upr L\vU^K), 
\endalign 
$$ 
by 5.4 Proposition, 2.6 Proposition respectively. Therefore 
$$ 
\align 
\Bbb A_F(\vT,\vU) &= \Psi_\vT\big(\Delta_F(\upr L\vT,\upr L\vU)\big) = \Psi_\vT\big(e^{-1}\Delta_K(\upr L\vT^K,\upr L\vU^K) \big) \\ &= \Psi_\vT
(e^{-1} \Psi_{\vT^K}^{-1}\big(\Bbb A_K(\vT^K,\vU^K)\big). 
\endalign 
$$ 
We write $y = \Bbb A_F(\vT,\vU)$, to get 
$$ 
y = \Psi_\vT(e^{-1}\Psi_{\vT^K}^{-1}(ey)). 
\tag 7.1.2  
$$ 
The Density Lemma of 5.2 says that the set of values $y = \Bbb A_F(\vT,\vU)$, $\vU\in \Scr E(F)$, is dense on the positive real axis, so (7.1.2) holds for all $y>0$. Writing $z = \Psi_{\vT^K}^{-1}(ey)$, we get $e^{-1}\Psi_{\vT^K}(z) = \Psi_\vT(e^{-1}z)$, as required. \qed 
\enddemo 
\remark{Remark} 
Given $\vT \in \Scr E(F)$, the definitions in \cite{2} (or see \cite{7} 6.3) give a finite tame extension $K/F$ for which $\vT$ has a totally wild $K/F$-lift. The proposition therefore reduces the problem of computing $\Psi_\vT$ to the case where $\vT$ is {\it totally wild.} 
\par 
When $\vT$ is totally wild and $K/F$ is tamely ramified, there is a simple relation (4.6) connecting $\vF_{\vT}$ and $\vF_{\vT^K}$. Likewise for $\vS_{\upr L\vT}$ and $\vS_{\upr L\vT^K}$ (3.2). However, for general $\vT$, the relations between $\vF_{\vT}$ and $\vF_{\vT^K}$, and between $\vS_{\upr L\vT}$ and $\vS_{\upr L\vT^K}$, are rather intricate. The symmetry indicated by the proposition can be viewed as a refinement of the Tame Parameter Theorem of \cite{7} 6.3. 
\endremark 
\subhead 
7.2 
\endsubhead 
Recall that $\sigma\in \wW F$ is {\it totally wild\/} if $\sigma|_{\scr P_F}$ is irreducible. Equivalently, the orbit $[\sigma,0]^+ \in \scr W_F \backslash \wP F$ has exactly one element. Write $\wwr F$ for the set of totally wild classes in $\wW F$. In particular, any $\sigma\in \wwr F$ has dimension $p^r$, for some $r\ge0$. 
\proclaim{Lemma} 
A representation $\sigma\in \wW F$ is totally wild if and only if $\sigma = \upr L\pi$, for $\pi \in \wG F$ such that $\roman{gr}(\pi) = \deg\vT_\pi$ and $\vT_\pi$  is totally wild. 
\endproclaim 
\demo{Proof} 
This follows from \cite{7} 6.3. \qed 
\enddemo 
\subhead 
7.3 
\endsubhead 
Totally wild representations of $\scr W_F$ exhibit simple ultrametric behaviour with respect to twisting by characters.  
\proclaim{Proposition} 
Let $\sigma \in \wwr F$ and let $c$ be a positive integer. If $\chi$ is a character of $\scr W_F$ of conductor $c$, then $\Delta(\sigma,\chi\otimes\sigma) \le c$. If $\vS_\sigma'$ is continuous at $c$, then $\vD(\sigma,\chi\otimes\sigma) = c$. 
\endproclaim 
\demo{Proof} 
Suppose $c>\roman{sl}(\sigma)$. The definition of $\vS_\sigma$ (3.1.2) shows that $\vS'_\sigma$ is continuous at $c$. Also $\roman{sl} (\chi\otimes \sigma) = c > \roman{sl}(\sigma)$, whence $\Delta(\sigma,\chi\otimes\sigma) = c$. We assume, therefore, that $c\le \roman{sl}(\sigma)$. The representations $\sigma$, $\chi\otimes \sigma$ are $\scr R_F^+(c)$-isomorphic, so $\Delta(\sigma,\chi\otimes \sigma) \le c$. The distance $\Delta(\sigma,\chi\otimes \sigma)$ is strictly less than $c$ if and only if $\chi|_{\scr R_F(c)}$ occurs in $\check\sigma \otimes \sigma|_{\scr R_F(c)}$. Suppose this condition holds. Since $\chi$ is trivial on $\scr R_F^+(c)$, the definition now shows that $\vS'_\sigma$ is discontinuous at $c$. \qed 
\enddemo 
\subhead 
7.4 
\endsubhead 
We recall how the set $\Scr E(F)$ carries a canonical action of the group of characters of $F^\times$. 
\par 
Let $\vT \in \Scr E(F)$, and let $\chi$ be a character of $F^\times$. If $\deg\vT = 1$, then $\vT$ is the endo-class of a character $\theta$ of $U^1_F$ and $\chi\vT$ is the endo-class of the (possibly trivial) character $\theta\chi\,\big|\,{U^1_F}$. Assume that $\deg\vT >1$ and choose a realization $\theta\in \scr C(\frak a,\beta)$ of $\vT$, relative to a simple stratum $[\frak a,m,0,\beta]$ in a matrix algebra $\End FV$. Define a character $\chi\theta$ of $H^1(\beta,\frak a)$ by 
$$
\chi\theta(h) = \chi(\det h)\,\theta(h),\quad h\in H^1(\beta,\frak a). 
$$ 
\proclaim{Lemma} 
Let $k = \sw(\chi) \ge 1$ and let $c\in F^\times$ satisfy $\chi(1{+}x) = \psi_F(cx)$, for $2\ups_F(x) > k$. If $m' = \roman{max}\,\{m,nk\}$, the quadruple $[\frak a,m',0,\beta{+}c]$ is a simple stratum and $\chi\theta\in \scr C(\frak a,\beta{+}c)$. 
\endproclaim 
\demo{Proof} 
See \cite{10} Appendix. \qed 
\enddemo 
Denote by $\chi\vT$ the endo-class of $\chi\theta$. If $\pi\in \Ao nF$ and $\vT = \vT_\pi$, then $\chi\vT$ is the endo-class $\vT_{\chi\pi}$ of the representation $\chi\pi:g\mapsto \chi(\det g)\pi(g)$, $g\in \GL nF$. 
\par 
The lemma shows that if $\vT$ is totally wild then so is $\chi\vT$, for any $\chi$. In a more general setting, the following is a direct consequence of the definitions (4.4.1), (3.1.2), on noting that $\upr L(\chi\vT) = \chi\otimes \upr L\vT$ (in the obvious notation). 
\proclaim{Proposition} 
Let $\vT \in \Scr E(F)$. If $\chi$ is a character of $F^\times$, then $\vF_{\chi\vT} = \vF_\vT$ and $\vS_{\upr L(\chi\vT)} = \vS_{\upr L\vT}$. Consequently, $\Psi_{\chi\vT} = \Psi_\vT$. 
\endproclaim 
\subhead 
7.5 
\endsubhead 
Our main result gives a procedure for calculating the Herbrand function $\Psi_\vT$ of any $\vT \in \Scr E(F)$. As noted in 7.1, it is enough to treat the case where $\vT$ is totally wild. 
\par 
If $\vT \in \Scr E(F)$ is totally wild and if $K/F$ is a finite tame extension, let $\vT^K\in \Scr E(K)$ be the unique $K/F$-lift of $\vT$. Denote by $\Bbb A_K$ the canonical ultrametric on $\Scr E(K)$. 
\proclaim{Interpolation Theorem} 
Let $\vT \in \Scr E(F)$ be totally wild. Let $D$ be a finite set of positive real numbers, containing all discontinuities of the functions $\Psi'_\vT$, $\vS_{\upr L\vT}'$. The function $\Psi_\vT$ has the following properties. 
\roster 
\item 
It is continuous, strictly increasing and piece-wise linear. 
\item 
Its derivative $\Psi_\vT'$ is continuous outside of $D$.  
\item 
If $K/F$ is a finite tame extension with $e = e(K|F)$ and if $\chi$ is a character of $K^\times$ satisfying $e^{-1} \sw(\chi) \notin D$, then 
$$ 
\Bbb A_K(\vT^K,\chi\vT^K) = e\Psi_\vT(e^{-1}\sw(\chi)). 
\tag 7.5.1 
$$ 
\endroster 
These properties determine $\Psi_\vT$ uniquely. 
\endproclaim 
\demo{Proof} 
The function $\Psi_\vT$ certainly satisfies (1) by 6.2 Proposition and (2) by definition of $D$. Condition (3) determines $\Psi_\vT(x)$ at a set of points $x$ dense in the positive real axis. Since $\Psi_\vT$ is continuous, it is thereby determined completely. 
\par 
We have to show that $\Psi_\vT$ has property (3). Let $\deg\vT = p^r$ and let $\pi \in \Ao{p^r}F$ to satisfy $\vT_\pi = \vT$. Set $\sigma = \upr L\pi \in \wwr F$. If $K/F$ is a finite tame extension, set $\sigma^K = \sigma|_{\scr W_K}$ and define $\pi^K\in \Ao{p^r}K$ by $\upr L\pi^K = \sigma^K$. In particular, $\vT_{\pi^K} = \vT^K$ \cite{7} 6.2 Proposition. 
\par 
Let $\chi$ be a character of $K^\times$, of conductor $k\ge 1$, such that $\vS_\sigma'$ is continuous at $k/e$. By 3.2 Proposition $\vS_{\sigma^K}(x) = e\vS_\sigma(x/e)$, so $\vS_{\sigma^K}'$ is continuous at $k$. By (5.4.1),  $\vs(\check\pi^K\times \chi\pi^K) = \vF_{\vT^K}(\Bbb A_K(\vT^K,\chi\vT^K))$. By 7.3 Proposition, 
$$ 
\vs(\check\sigma^K\otimes \chi\otimes\sigma^K) = \vS_{\sigma^K}(\Delta_K(\sigma^K,\chi\otimes\sigma^K)) = \vS_{\sigma^K}(k),  
$$ 
where $\Delta_K$ is the canonical pairing on $\wW K$. By 4.6 Proposition, $\vF_{\vT^K}(x) = e\vF_\vT(x/e)$. Altogether 
$$ 
\align 
\vs(\check\pi^K\times \chi\pi^K) &= \vF_{\vT^K}\big(\Bbb A_K(\chi\vT^K,\vT^K)\big)  = e\,\vF_\vT\big(e^{-1} \Bbb A_K(\chi\vT^K,\vT^K)\big) \\ 
&= \vs(\check\sigma^K\otimes \chi\otimes\sigma^K) = \vS_{\sigma^K}(k) \\ 
&= e\vS_\sigma(k/e), 
\endalign 
$$ 
whence $\Bbb A_K(\chi\vT^K,\vT^K) = e\,\vF_\vT^{-1}\circ \vS_\sigma(k/e) = e\Psi_\vT(k/e)$, as required. \qed 
\enddemo 
The proof shows that (7.5.1) holds provided only that $\vS_{\upr L\vT}'$ is continuous at $\sw(\chi)/e$. We have no example where $\Psi_\vT'$ is continuous at a discontinuity of $\vS_{\upr L\vT}'$. However, in practice, there is no advantage in minimizing the finite set $D$ of ``exceptional'' points. 
\subhead 
7.6 
\endsubhead 
We describe the function $\Psi_\vT$, for $\vT \in \Scr E(F)$ totally wild, on part of its range. We have already noted 
(6.2) that $\Psi_\vT(0) = 0$ and that $\Psi_\vT(x) = x$ for $x> m_\vT$. 
\proclaim{Proposition} 
Let $\vT \in \Scr E(F)$ be totally wild of degree $p^r$, $r\ge0$, and write $m_\vT = ap^{t-r}$, for integers $a$, $t$ with $a \not\equiv 0 \pmod 
p$ and $0\le t < r$. 
\roster 
\item There exists $\eps > 0$ such that $\Psi_\vT'(x) = p^{-r}$, for $0< x < \eps$. 
\item There exists $\delta > 0$ such that 
$$ 
\Psi'_\vT(x) = p^{r-t}, \quad m_\vT{-}\delta < x < m_\vT. 
$$ 
\endroster 
\endproclaim 
\demo{Proof} 
Part (1) follows from the definitions (4.4.1) and (3.2.1) on noting that, if $\sigma\in \wwr F$, there exists $\eps >0$ such that $\sigma$ is irreducible on $\scr R_F(\eps)$. 
\par
In part (2), write $m_\vT = p^{-r}m = p^{t-r}a$. The class $\vT$ has a realization $\theta \in \scr C(\frak a,\beta)$, for a simple stratum $[\frak 
a,m,0,\beta]$ in $\M{p^r}F$. We choose a simple stratum $[\frak a,m,m{-}1,\alpha]$ equivalent to $[\frak a,m,m{-}1,\beta]$. The degree 
$[F[\alpha]{:}F]$ is then $p^{r-t}$. Thus $\vF_\vT'(x) = p^{t-r}$ in a region $m_\vT{-}\delta < x < m_\vT$. 
\par 
Since $\roman{sl}(\sigma) = \vs(\sigma) = m_\vT$, the restriction of $\sigma$ to $\scr R_F(m_\vT)$ is a sum of characters of $\scr R_F(m_\vT)$ 
trivial on $\scr R_F^+(m_\vT)$. These are all conjugate under $\scr P_F$ and so, by 1.2 Lemma 2, they are all the same. Therefore, 
every irreducible component of $\check\sigma\otimes\sigma$ contains the trivial character of $\scr R_F(m_\vT)$. By 2.1 Proposition 1, $
\check\sigma\otimes\sigma$ is trivial on $\scr R_F(m_\vT{-}\delta)$, for some $\delta > 0$. In that region, $\vS'_\sigma$ has value $1$, 
whence the result follows. \qed 
\enddemo 
\remark{Remark} 
In the context of the proposition, one can have $m_\vT = ap^{t-r}$, with $a\not\equiv 0 \pmod p$ and $t\ge r$. This, however, is equivalent to the existence of a character $\chi$ of $F^\times$ such that $m_{\chi\vT} < m_\vT$. In light of 7.4 Proposition, nothing is lost by excluding this case. 
\endremark 
\example{7.7 Example} 
Say that $\vT \in \Scr E(F)$ is {\it essentially tame} if, for some finite, tamely ramified extension $K/F$, $\vT$ has a $K/F$-lift of degree $1$: equivalently, $e_\vT$ is relatively prime to $p$. 
\proclaim{Corollary} 
An endo-class $\vT \in \Scr E(F)$ satisfies $\Psi_\vT(x) = x$, $x\ge 0$, if and only if $\vT$ is essentially tame. 
\endproclaim 
\demo{Proof} 
7.6 Proposition shows that a totally wild endo-class $\vX \in \Scr E(F)$ has the property $\Psi_\vX(x) = x$ if and only if $\deg\vX = 1$. Conversely, take $\vT \in \Scr E(F)$. the definition in \cite{2} shows that there exists a finite, tame, Galois extension $K/F$ such that $\vT$ has a totally wild $K/F$-lift $\vT^K$. By 7.1 Proposition, $\Psi_\vT(x) = x$ if and only if $\Psi_{\vT^K}(x) = x$, whence the result follows. \qed 
\enddemo 
\endexample 
\head\Rm 
8. The decomposition function 
\endhead 
We analyze some features of the decomposition function $\vS_\sigma$, $\sigma \in \wW F$, taking the view that $\vS_\sigma$ has been given somehow, without prior knowledge of $\sigma$. 
\subhead 
8.1 
\endsubhead 
We examine the discontinuities of the derivative $\vS_\sigma'$, using only group-theoretic methods. 
\proclaim{Proposition} 
Let $\sigma \in \wW F$, let $\eps>0$ and let $\sigma_\eps$ be an irreducible component of $\sigma\,\big|\,\scr R_F(\eps)$. Let $\Gamma_\eps$ be the group of characters of $\scr R_F(\eps)/\scr R_F^+(\eps)$. The following are equivalent. 
\roster 
\item The function $\vS'_\sigma$ is continuous at $\eps$. 
\item The representation $\chi\otimes\sigma_\eps$ is not $\scr W_F$-conjugate to $\sigma_\eps$, for any $\chi\in \Gamma_\eps$, $\chi\neq1$. 
\endroster 
\endproclaim 
\demo{Proof} 
An  exercise in elementary representation theory yields: 
\proclaim{Lemma} 
Suppose that the representation $\sigma_\eps\,\big|\,{\scr R_F^+(\eps)} = \sigma_\eps^+$ is irreducible. The map $\chi\mapsto \chi\otimes \sigma_\eps$ is a bijection between the group $\Gamma_\eps$ and the set of isomorphism classes of irreducible smooth representations of $\scr R_F(\eps)$ that contain $\sigma_\eps^+$. 
\endproclaim 
We prove the proposition. For $\delta >0$, define 
$$ 
\align 
d_\delta &= \dim\Hom{\scr R_F(\delta)}1{\check\sigma\otimes\sigma}, \\ 
d_\delta^+ &= \dim\Hom{\scr R^+_F(\delta)}1{\check\sigma\otimes\sigma}. 
\endalign 
$$ 
The step function $\vS'_\sigma$ is continuous at a point $\delta >0$ if and only if it is constant on a neighbourhood of $\delta$. This is equivalent to the condition $d_\delta = d^+_\delta$. 
\par 
Let $m_\delta$ be the multiplicity of $\sigma_\delta$ in $\sigma\,\big|\,{\scr R_F(\delta)}$ and $l_\delta$ the number of $\scr W_F$-conjugates of $\sigma_\delta$. Define $m_\delta^+$ and $l_\delta^+$ analogously, relative to an irreducible component $\sigma_\delta^+$ of $\sigma_\delta\,\big|\,{\scr R_F^+(\delta)}$. Thus $d_\delta = l_\delta m_\delta^2$ and $d^+_\delta = l^+_\delta{m^+_\delta}^2$. Moreover, $l_\delta m_\delta$ and $l^+_\delta m_\delta^+$ are the Jordan-H\"older lengths of the restrictions $\sigma\,\big|\,{\scr R_F(\delta)}$ and $\sigma\,\big|\,{\scr R_F^+(\delta)}$ respectively. 
\par 
Suppose that condition (2) holds. In particular, $\sigma_\eps \not\cong \sigma_\eps \otimes \chi$, for any character $\chi\in \Gamma_\eps$. This implies that $\sigma_\eps\,\big|\,\scr R_F^+(\eps)$ is irreducible. Writing $\sigma^+_\eps= \sigma_\eps\,\big|\,{\scr R_F^+(\eps)}$, the lemma says that $\sigma_\eps$ is the unique irreducible component of $\sigma\,\big|\,{\scr R_F(\eps)}$ containing $\sigma_\eps^+$. Thus $m_\eps = m_\eps^+$. The representations $\sigma_\eps$, $\sigma^+_\eps$ have the same $\scr W_F$-isotropy, so $l^+_\eps = l_\eps$. Therefore $d_\eps = d^+_\eps$ and $\vS_\sigma'$ is continuous at $\eps$. 
\par 
Suppose now that (2) fails. If $\sigma_\eps\,\big|\,\scr R_F^+(\eps)$ is reducible, certainly $\vS'_\sigma$ cannot be continuous at $\eps$. We therefore assume the contrary and let $c$ be the number of $\chi\in \Gamma_\eps$ such that $\chi\otimes \sigma_\eps$ is $\scr W_F$-conjugate to $\sigma_\eps$. Thus $c>1$ and, by the lemma, $l_\eps = c\,l_\eps^+$. Correspondingly, $m_\eps^+ = cm_\eps$, so $d_\eps^+ = cd_\eps > d_\eps$ and $\vS'_\sigma$ is not locally constant at $\eps$. \qed 
\enddemo 
\remark{Remark} 
We draw attention to one step in the preceding proof: if the conditions of the proposition are satisfied, then $\sigma_\eps\,\big|\,\scr R_F^+(\eps)$ is {\it irreducible.} 
\endremark 
\subhead 
8.2 
\endsubhead 
To prepare for the main result, we need some ideas from Galois theory. 
\par 
Let $\sigma \in \wW F$, and assume $\dim\sigma >1$. Define $\bar\sigma$ to be the {\it projective\/} representation defined  by $\sigma$: that is, if $\dim\sigma = n$, then $\bar\sigma$ is the composition of $\sigma$ with the canonical map $\GL n{\Bbb C} \to \roman{PGL}_n(\Bbb C)$. The image of $\bar\sigma$ is finite and $\roman{Ker}\,\bar\sigma$ is of the form $\scr W_E$, for a finite Galois extension $E/F$. We call $E/F$ the {\it pro-kernel field\/} of $\sigma$. Let $T/F$ be the maximal tamely ramified sub-extension of $E/F$: we call $T/F$ the {\it tame kernel field\/} of $\sigma$. 
\definition{Definition}  
Let $\sigma\in \wwr F$. 
\roster 
\item 
Define $D(\sigma)$ as the group of characters $\chi$ of $\scr W_F$ such that $\chi\otimes \sigma \cong \sigma$. 
\item 
Write $\sigma_0^+ = \sigma\,\big|\,{\scr P_F} \in \wP F$. Define $D_0(\sigma)$ as the group of characters $\phi$ of $\scr P_F$ such that $\phi\otimes \sigma_0^+ \cong \sigma_0^+$. 
\endroster 
\enddefinition   
Restriction of characters gives a canonical homomorphism $D(\sigma) \to D_0(\sigma)$. A character $\phi$ of $\scr P_F$ lies in $D_0(\sigma)$ if and only if it is a component of $\check\sigma_0^+\otimes \sigma_0^+$, whence $|D_0(\sigma)| \le (\dim\sigma)^2$. On the other hand, $\sigma_0^+$ is effectively an irreducible representation of a finite $p$-group, of dimension $>1$. Consequently, the group $D_0(\sigma)$ is not trivial. 
\par 
The representation $\sigma_0^+$ is stable under conjugation by $\scr W_F$, so $\scr W_F$ acts on $D_0(\sigma)$, with $\scr P_F$ acting trivially. The $\scr W_F$-stabilizer of a character $\phi \in D_0(\sigma)$ is thus of the form $\scr W_{T_\phi}$, for a finite tame extension $T_\phi/F$. The kernel of the canonical map $\scr W_F \to \roman{Aut}\,D_0(\sigma)$ is therefore 
$$ 
\scr W_{T_I} = \bigcap_{\phi\in D_0(\sigma)} \scr W_{T_\phi},  
$$ 
where $T_I/F$ is a finite, tamely ramified, Galois extension. We call $T_I/F$ the {\it imprimitivity field\/} of $\sigma$. 
\par 
If $K/F$ is a finite tame extension, the representation $\sigma^K = \sigma\,\big|\,{\scr W_K}$ is irreducible and lies in $\wwr K$. It agrees with $\sigma$ on $\scr P_K = \scr P_F$ so $D_0(\sigma^K) = D_0(\sigma)$.  
\proclaim{Proposition} 
If $\sigma\in \wwr F$ has tame kernel field $T/F$ and imprimitivity field $T_I/F$, then $T_I \subset T$. The canonical map $D(\sigma^{T_I}) \to D_0(\sigma)$ is an isomorphism. 
\endproclaim 
\demo{Proof} 
We first note: 
\proclaim{Lemma 1} 
If $\zeta\in D(\sigma)$ is tamely ramified then $\zeta = 1$. 
\endproclaim 
\demo{Proof} 
The kernel of $\zeta$ is $\scr W_K$, for a finite, cyclic, tame extension $K/F$. The relation $\zeta\otimes\sigma \cong \sigma$ implies that $\sigma$ is reducible on $\scr W_K$. Since $\scr P_F \subset \scr W_K$, it is also reducible on $\scr P_F$, contrary to hypothesis. \qed 
\enddemo 
\proclaim{Lemma 2} 
Let $K/F$ be a finite tame extension. The restriction map $D(\sigma^K) \to D_0(\sigma)$ is an isomorphism of $D(\sigma^K)$ with the group $D_0(\sigma)^{\scr W_K}$ of $\scr W_K$-fixed points in $D_0(\sigma)$. 
\endproclaim 
\demo{Proof} 
Lemma 1 implies that the map $D(\sigma^K) \to D_0(\sigma)$ is injective. Its image is clearly contained in $D_0(\sigma)^{\scr W_K}$. Let $\zeta\in D_0(\sigma)^{\scr W_K}$. Thus $\zeta$ admits extension to a character $\tilde\zeta$ of $\scr W_K$ \cite{7} 1.3 Proposition. The representations $\sigma^K$, $\tilde\zeta\otimes\sigma^K$ agree on $\scr P_K$ so ({\it loc\. cit\.}) there is a tame character $\chi$ of $\scr W_K$ such that $\chi\tilde\zeta\otimes\sigma^K \cong \sigma^K$. Therefore $\chi\tilde\zeta\in D(\sigma^K)$, as required. \qed 
\enddemo 
By the definition of $T$, we have $\scr W_T = \scr P_F\scr W_E$. A character $\zeta\in D_0(\sigma)$ is effectively a character of $\scr P_F\scr W_E/\scr W_E$, and hence a character of $\scr W_T$. In particular, $\scr W_T$ fixes $\zeta$, whence $\scr W_T \subset \scr W_{T_\zeta}$. Therefore $\scr W_T \subset \scr W_{T_I}$, or $T\supset T_I$, as required to complete the proof of the proposition. \qed 
\enddemo 
\remark{Remark} 
There are examples of representations $\sigma\in \wwr F$ such that $T_I \neq T \neq E$. 
\endremark 
\subhead 
8.3 
\endsubhead 
Let $\sigma\in \wwr F$. Say that $\sigma$ is {\it absolutely wild\/} if its tame kernel field is $F$. That is, if $E$ is the pro-kernel field of $\sigma$, then $E/F$ is totally wildly ramified. 
\proclaim{Theorem} 
Let $\sigma\in \wwr F$ be absolutely wild of dimension $p^r$, $r\ge1$. If $a>0$ is the least discontinuity of $\vS_\sigma'$ then $a$ is an integer and 
$$ 
a = \roman{min}\{\sw(\chi): \chi\in D(\sigma), \chi\neq1\}. 
$$ 
\endproclaim 
\demo{Proof} 
Nothing is changed by tensoring $\sigma$ with a tame character of $\scr W_F$. We may therefore assume that $\sigma$ is a representation of $\Gal {\widetilde E}F$, where $\widetilde E/F$ is totally wildly ramified. 
\par 
The group $D(\sigma) \cong D_0(\sigma)$ is non-trivial. We accordingly define 
$$ 
c = \roman{min}\,\{\sw(\chi): \chi\in D(\sigma), \chi\neq1\}. 
$$ 
Suppose first that $c<a$: in particular, $\vS'_\sigma$ is continuous at $c$. Any character $\phi \in D(\sigma)$ occurs as an irreducible component of $\check\sigma\otimes\sigma$, so the definitions in 5.1 imply that $\vS'_\sigma$ is discontinuous at $c$. This contradiction implies $c\ge a$. 
\par 
If $0<\eps < a$, condition (2) of 8.1 Proposition holds at $\eps$, so $\sigma_\eps = \sigma\,\big|\,\scr R_F(\eps)$ is irreducible (8.1 Remark). It follows that $\sigma_a = \sigma\,\big|\,\scr R_F(a)$ is also irreducible. Since $\vS_\sigma'$ is discontinuous at $a$, there is a non-trivial character $\chi$ of $\scr R_F(a)/\scr R_F^+(a)$ such that $\sigma_a\otimes \chi$ is $\scr W_F$-conjugate to $\sigma_a$, say $\sigma_a^g \cong \sigma_a\otimes \chi$ for some $g\in \scr W_F$. However, $\sigma_a = \sigma\,\big|\,\scr R_F(a)$ and surely $\sigma^g \cong \sigma$. Thus $\sigma_a \cong \sigma_a\otimes \chi$, whence $\sigma_a$ is reducible on $\scr R_F^+(a)$. Consequently, $\sigma$ is reducible on $\scr R_F^+(a)$. As $\sigma$ is effectively a representation of a finite $p$-group, it is induced from a representation of an open normal subgroup of $\scr W_F$, of index $p$ and containing $\scr R_F^+(a)$. That is, there is a non-trivial character $\phi$ of $\scr W_F$, trivial on $\scr R_F^+(a)$, such that $\sigma\otimes \phi \cong \sigma$. Therefore $c\le \sw(\phi) \le a$, giving $c=a$, as required. \qed 
\enddemo 
The proof of the theorem relies on $\sigma$ being absolutely wild, but the result extends to the general case of $\sigma\in \wwr F$. 
\proclaim{Corollary} 
Let $\sigma \in \wwr F$ have dimension $p^r$, $r\ge1$. Let $T_I/F$ be the imprimitivity field of $\sigma$ and set $e = e(T_I|F)$. The least discontinuity $a$ of $\vS_\sigma'$ is given by  
$$
a = \roman{min}\,\{\sw(\chi)/e:\chi\in D(\sigma^{T_I}),\chi\neq 1\}. 
$$ 
In particular, $a$ is $p$-integral. 
\endproclaim 
\demo{Proof} 
We apply the theorem to the absolutely wild representation $\sigma^T\in \wwr T$, where $T/F$ is the tame kernel field of $\sigma$. If $c$ is the least discontinuity of $\vS_\sigma'$, then $e(T|F)c$ is that of $\vS'_{\sigma^T}$. If $\phi\in D(\sigma^T)$, then $\phi = \chi|_{\scr W_T}$, for a unique $\chi\in D(\sigma^{T_I})$ (8.2 Proposition), and $\sw(\phi) = e(T|T_I)\,\sw(\chi)$. \qed 
\enddemo 
\subhead 
8.4 
\endsubhead 
Similar techniques lead to an attractive result, which we leave as an exercise.  
\proclaim{Proposition} 
Let $\sigma\in \wW F$ be absolutely wild of dimension $p^r$, $r\ge1$, with pro-kernel field $E/F$. If $\vS_\sigma'$ has only one discontinuity, then $\Gal EF$ is elementary abelian of order $p^{2r}$. 
\endproclaim 
\subhead 
8.5 
\endsubhead 
Consider, as an example, the case where $\deg\vT = p$, $\vT \in \Scr E(F)$. Thus $\vT$ is either essentially tame or  totally wild. The first case is covered by 7.7, so assume $\vT$ totally wild. Write $m_\vT = m/p$. Twisting with a character of $F^\times$ changes nothing, so we assume $m\not\equiv 0 \pmod p$. 
\par  
Directly from (4.4.1) and 4.1 Proposition, we have $\vF_\vT(0) = m(p{-}1)/p^2$ and $\vF_\vT'(x) = p^{-1}$, $0<x<m/p$. On the other side, $\vS_{\upr L\vT}(0) = \vF_\vT(0)$. The only possible values for $\vS_{\upr L\vT}'$ are $p^{-2}$, $p^{-1}$ and $1$. By 7.6 Proposition, the first and third certainly occur. If only they occur, we are in the case of 8.4. The unique discontinuity of $\vS'_{\upr L\vT}(x)$ occurs at $x=m/(p{+}1)$. Otherwise, there are two discontinuities, the first of which is interpreted by 8.3 Corollary and that determines the second one. This case is analyzed in detail by M\oe glin \cite{20}. 
\subhead 
8.6 
\endsubhead  
By way of a contrast, take $p=2$ and assume that $F$ contains a primitive cube root of unity. We are aware of an example of a totally wild endo-class $\vT\in \Scr E(F)$ with $\deg\vT = 4$, $m_\vT = 1/2$, having the following properties. First, $\vS'_{\upr L\vT}$ has a unique discontinuity at $1/3$, while $\vF'_\vT$ has a discontinuity at $1/4$ and $\vF_\vT(0) = 5/16$. So, in the interesting range $0<x<m_\vT = 1/2$, the function $\Psi'_\vT$ has discontinuities at $1/3$ and $3/8$. In all, 
$$ 
\Psi_\vT'(x) = \left\{\,\alignedat3 &\tfrac14, &\quad &0<x<\tfrac13, \\ &4, &\quad &\tfrac13<x<\tfrac38, \\ &2, &\quad &\tfrac38 < x< \tfrac12. \endalignedat \right. 
$$ 
(One takes a tensor product $\sigma = \sigma_1\otimes\sigma_2$, where the $\sigma_i \in \wwr F$ are carefully chosen of dimension $2$ and Swan conductor $1$. The endo-class $\vT$ is given by $\upr L\vT = [\sigma;0]^+$.) 
 

\Refs 
\ref\no 1 
\by C.J. Bushnell 
\paper  Effective local Langlands correspondence 
\inbook Automorphic forms and Galois representations vol. 1 \eds F. Diamond, P. L Kassei, Minhyong Kim \bookinfo London Math. Soc. Lecture Notes \vol 414 \publ Cambridge University Press \yr 2014 \pages 102--134 
\endref 
\ref\no2 
\by C.J. Bushnell and G. Henniart 
\paper Local tame lifting for $\roman{GL}(n)$ I: simple characters 
\jour Publ. Math. IHES \vol 83 \yr 1996 \pages 105--233 
\endref 
\ref\no3
\bysame 
\paper Local tame lifting for $\roman{GL}(n)$ II: wildly ramified supercuspidals 
\jour Ast\'erisque \vol 254 \yr 1999 \pages vi+105 
\endref 
\ref\no4
\bysame 
\paper Local tame lifting for $\roman{GL}(n)$ IV: simple characters and base change 
\jour Proc. London Math. Soc. (3) \vol 87 \yr 2003 \pages 337--362 
\endref 
\ref\no5 
\bysame 
\book The local Langlands Conjecture for $\roman{GL}(2)$ 
\bookinfo Grundlehren der mathematischen Wissenschaften {\bf 335} \publ Springer \yr 2006 
\endref 
\ref\no6
\bysame 
\paper Intertwining of simple characters in $\roman{GL}(n)$ \jour Int. Math. Res. Not. IMRN \yr 2013 \vol 17 \pages 3977--3987 
\endref  
\ref\no7 
\bysame 
\paper To an effective local Langlands Correspondence \jour Memoirs Amer. Math. Soc.  \issue 1087 \vol {231} \yr 2014  \pages v+88 
\endref 
\ref\no8 
\by C.J. Bushnell, G. Henniart and P.C. Kutzko 
\paper Local Rankin-Selberg convolutions for $\roman{GL}_n$: Explicit conductor formula 
\jour J. Amer. Math. Soc. \vol 11 \yr 1998 \pages 703--730 
\endref 
\ref\no9 
\by C.J. Bushnell and P.C. Kutzko 
\book The admissible dual of $GL(N)$ via compact open subgroups 
\bookinfo Annals of Math. Studies {\bf 129} \publ Princeton University Press \yr 1993 
\endref 
\ref\no10 
\bysame 
\paper The admissible dual of $SL(N)$ II \jour  Proc. London Math. Soc. (3) \vol 68 \yr 1994 \pages 317--379 
\endref 
\ref\no11  
\bysame 
\paper  Simple types in $GL(N)$: computing conjugacy classes \inbook 
Representation theory and analysis on homogeneous spaces \ed S. Gindikin et al. \bookinfo Contemp. Math. \vol 177 \yr 1994 \pages 107--135 
\endref 
\ref\no12 
\bysame 
\paper Smooth representations of $p$-adic reductive groups; Structure theory via types \jour Proc. London Math. Soc.  (3) \vol 77 \yr 1998 \pages 582--634 
\endref 
\ref\no13 
\bysame 
\paper Semisimple types for $GL(N)$ \jour Compositio Math. \vol 119 \yr 1999 \pages 53--97 
\endref \ref\no14 
\by M. Harris and R. Taylor 
\book On the geometry and cohomology of some simple Shimura varieties 
\bookinfo Annals of Math. Studies {\bf 151} \publ Princeton University Press \yr 2001 
\endref 
\ref\no15  
\by V. Heiermann 
\paper Sur l'espace des repr\'esentations irr\'eductibles du groupe de Galois d'un corps local \jour C. R. Acad. Sci. Paris S\'er. I Math. \vol 323 \issue 6 \yr 1996 \pages 571--576 
\endref 
\ref\no16 
\by  G. Henniart 
\paper Repr\'esentations du groupe de Weil d'un corps local  \jour L'Ens. Math. S\'er II \vol 26 \yr 1980 \pages 155-172 
\endref 
\ref\no17 
\bysame 
\paper Une preuve simple des conjectures locales de Langlands pour $\roman{GL}_n$ sur un corps $p$-adique 
\jour Invent. Math. \vol 139 \yr 2000 \pages 439--455
\endref 
\ref\no18 
\by H. Jacquet, I. Piatetski-Shapiro and J. Shalika 
\paper Rankin-Selberg convolutions 
\jour Amer. J. Math. \vol 105 \yr 1983 \pages 367--483 
\endref 
\ref\no19
\by G. Laumon, M. Rapoport and U. Stuhler 
\paper $\Cal D$-elliptic sheaves and the Langlands correspondence 
\jour Invent. Math. \vol 113 \yr 1993 \pages 217--338 
\endref 
\ref\no20 
\by C. M\oe glin 
\paper Sur la correspondance de Langlands-Kazhdan 
\jour J. Math. Pures et Appl. (9) \vol 69 \yr 1990 \pages 175--226 
\endref 
\ref\no21 
\by P. Scholze 
\paper The local Langlands correspondence for $\roman{GL}_n$ over $p$-adic fields 
\jour Invent. Math. \vol 192 \yr 2013 \pages 663--715 
\endref 
\ref\no22 
\by J-P. Serre 
\book Corps locaux \publ Hermann \publaddr Paris \yr 1968 
\endref 
\ref\no23 
\by F. Shahidi  
\paper Fourier transforms of intertwining operators and Plancherel measures for $\roman{GL}(n)$ 
\jour Amer. J. Math. \vol 106 \yr 1984 \pages 67--111
\endref  
\endRefs
\enddocument